\documentclass[titlepage,11pt]{article}
\usepackage{latexsym}
\usepackage{amsfonts}
\usepackage{amsmath}
\usepackage{comment}
\newcommand\blackslug{\hbox{\hskip 1pt \vrule width 4pt height 8pt depth 1.5pt
        \hskip 1pt}}
\newcommand\bbox{\hfill \quad \blackslug \bigbreak}

\def\DD{\hbox{-}}
\def\CC{\hbox{-}\cdots\hbox{-}}
\def\LL{,\ldots,}

\def\MM{\mathcal{L}}
\def\dom{\operatorname{dom}}
\def\edom{\operatorname{edom}}

\title{Some results and problems on tournament structure}
\author{Tung Nguyen\thanks{Supported by AFOSR grant  FA9550-22-1-0234 and by NSF grant  DMS-2154169.}\\
Princeton University,\\ Princeton, NJ 08544, USA
\and
Alex Scott\thanks{Research supported by EPSRC grant EP/X013642/1.}\\
Mathematical Institute,\\ University of Oxford,
\\
Oxford OX2 6GG, UK
\and
Paul Seymour\footnotemark[1]\\
Princeton University,\\ Princeton, NJ 08544, USA}

\date{February 24, 2023; revised \today}

\newtheorem{thm}{}[section]

\newcommand{\Proof}{\noindent{\bf Proof.}\ \ }

\begin{document}
\maketitle
\begin{abstract}
This paper is a survey of results and problems related to the following question: is it true that if $G$ is a tournament with
sufficiently large chromatic number, then $G$ has two vertex-disjoint subtournaments $A,B$, both with large chromatic number,
such that all edges between them are directed from $A$ to $B$? 
We describe what we know about this question, and report some progress on several other related questions, on tournament colouring
and domination.
\end{abstract}

\section{Introduction}

There is an open question of El-Zahar and Erd\H{o}s~\cite{elzahar}, the following:
\begin{thm}\label{erdosconj}
{\bf Problem: }Is the following true? For all integers $t, c\ge 1$,
there exists $d\ge 1$, such that if a graph $G$ satisfies $\chi(G)\ge d$ and has no clique with $t$ vertices,
then there are subsets $A,B\subseteq V(G)$ with $\chi(G[A]), \chi(G[B])\ge c$, such that there are no edges between $A$ and $B$.
\end{thm}
($\chi(G)$ denotes the chromatic number of a graph $G$.)
In this paper we consider an analogue of this for tournaments, and a number of related results and conjectures.
(Henceforth we use ``Conjecture'' for open questions: we do not intend to imply that we believe them to be true.)

The {\em chromatic number} $\chi(T)$ of a tournament $T$ is the minimum $k$ such that $V(T)$ is the union of $k$
subsets $A_1\LL A_k$ that each induce an acyclic subtournament. If $T$ is a tournament and $A\subseteq V(T)$, we denote by $T[A]$
the subtournament with vertex set $A$, and write $\chi(A)$ for $\chi(T[A])$. If $A,B$ are disjoint, 
and there is no edge with tail in $B$ and head in $A$,
we say $A$ is {\em complete to} $B$, and $B$ is {\em complete from} $A$, and write $A\Rightarrow B$, and call $(A,B)$ 
a {\em complete pair}. 

In~\cite{dense} we made the following conjecture, which implies \ref{erdosconj}:
\begin{thm}\label{tourconj}
{\bf Conjecture: }For all $c\ge 0$ there exists $d\ge 0$ such that if $T$ is a tournament with $\chi(T)\ge d$, there is a complete pair $(A,B)$ of $T$
such that $\chi(A),\chi(B)\ge c$.
\end{thm}
Very recently (since an earlier draft of this paper was submitted for publication) Klingelhoefer and Newman~\cite{klingel} proved the equivalence
of \ref{erdosconj} and \ref{tourconj}, and we have updated this paper accordingly. In particular, we give a different proof of the equivalence that is 
simpler and shorter than
that in~\cite{klingel}, although it uses some of the same ideas.

It seems to us that \ref{tourconj} is very strong, and may well be false, 
despite its equivalence with \ref{erdosconj}. Indeed, 
for several years it was an open conjecture even to show that a tournament with sufficiently large chromatic number
has a vertex whose out-neighbour set has large chromatic number, and the proof of this~\cite{harut} is highly non-trivial.

But there are some partial results in its favour.
First, a {\em cyclic triangle} is a tournament with three vertices,
each with out-degree one; and we use the same name for the vertex set of such a tournament. Thus, a tournament is transitive if and only if it contains no cyclic triangle. 
In defense of \ref{tourconj}, we will prove:
\begin{thm}\label{triout}
For all $c\ge 0$ there exists $d\ge 0$ such that if $T$ is a tournament with $\chi(T)\ge d$,  there is a complete pair
$(A,B)$ of $T$ such that $A$ is a cyclic triangle and $\chi(B)\ge c$.
\end{thm}

The {\em domination number} $\dom(T)$ of a tournament $T$ is the size of the smallest set $X$ of vertices such that every vertex in $V(T)\setminus X$
is adjacent from a vertex in $X$; it is always at most the chromatic number. A second result in defense of \ref{tourconj} is that it is true
for tournaments with sufficiently large domination number. More exactly, we will show:

\begin{thm}\label{bigdom0}
For every integer $c\ge 1$, there exists $d\ge 1$ such that if $T$ is a tournament with $\dom(T)\ge d$
then there is a complete pair $(A,B)$
such that $\chi(A),\chi(B)\ge c$.
\end{thm}

This area 
seems to be relatively unexplored, and yet full of interesting, significant, interconnected questions. Indeed, our attempts to decide \ref{tourconj}, 
while unavailing, led to progress on several other tournament problems, for instance:
\begin{itemize}
\item A {\em rebel} is a tournament $H$ such that all $H$-free tournaments have bounded domination number. (A tournament is 
{\em $H$-free}
if no subtournament is isomorphic to $H$.)
Until now there was only one tournament with chromatic number more than two that was known to be a rebel, 
and the proof that it was a rebel 
was difficult~\cite{tourdom}. We give a much more general and much simpler construction of rebels.
\item A {\em hero} is a tournament $H$ such that all $H$-free tournaments have bounded chromatic number. It was already known which 
tournament are heroes~\cite{heroes}, but our methods give a simpler proof of this result.
\item A {\em legend} is an ordered tournament that is contained in every ordered tournament with sufficiently large
domination number. (An {\em ordered tournament} is a pair $(T,\tau)$, where $T$ is a tournament
and $\tau$ is a numbering of its vertex set; and its {\em domination number} is the domination number of $T$.) We will find all legends.
\item Harutyunyan, Le, Thomass\'e, and Wu~\cite{harut} proposed two open questions about domination. We will show that one implies the other.
\item We give a new class of tournaments (``crossing tournaments'') with surprising properties, that provide
counterexamples to some of our wilder dreams.
\end{itemize} 
This paper is a survey of what we know about \ref{tourconj}, and its connections with some other questions.
It is organized as follows: 
\begin{itemize}
\item In the next section we survey a number of special cases and variants of \ref{tourconj}, true, false and open.
\item In section \ref{sec:diamonds} we prove a result about ``diamonds'' that will have several applications later.
\item The next three sections concern the equivalence of \ref{tourconj} with the problem of El-Zahar and Erd\H{o}s.
\item Sections \ref{sec:heroes} and \ref{sec:makingrebels} prove the result about legends, and make progress on understanding rebels.
\item Section \ref{sec:harutconj} proves an implication between two open conjectures of~\cite{harut}.
\item Sections \ref{sec:density} and \ref{sec:laws} give a proof of \ref{triout} and several related results.
\item Sections \ref{sec:1HH} and \ref{sec:St} concern whether excluding  two particular types of tournament forces a bound on
``local chromatic number''.
\item Finally, section \ref{sec:makingheroes} applies some results of earlier sections to deduce an older theorem about 
how to construct heroes.
\end{itemize}

\section{Complete pairs $(A,B)$ with $\chi(A)$ small}

The conjecture \ref{tourconj} says that if $\chi(T)$ is large then there is a complete pair $(A,B)$ with both $\chi(A),\chi(B)$ large.
Can we at least get a complete pair $(A,B)$ with $\chi(A),\chi(B)$ small but nonzero?  Or perhaps with one of them large?

It is not obvious even that there is a complete pair $(A,B)$ with $|A|=1$ and $\chi(B)$ large, and this was raised in~\cite{heroes}
as a conjecture. But Harutyunyan, Le, Thomass\'e, and Wu proved the following fundamental result, in a breakthrough paper~\cite{harut}:
\begin{thm}\label{bigandsmall}
For every integer $c\ge 1$, there exist integers $K,k\ge 1$ such that
every tournament $T$ with $\dom(T)\ge K$ contains a subtournament on at most $k$ vertices having chromatic number
at least $c$.
\end{thm}
This has several useful consequences, including a proof of the conjecture of~\cite{heroes}:
\begin{thm}\label{outnbrs}
For all integers $c\ge 0$ there exists $d\ge 0$ such that for every tournament $T$ with $\chi(T)\ge d$, there exists $v\in V(T)$ such that
$\chi(N^+(v))\ge c$.
\end{thm}
For a vertex $v$, we define $N^+[v]=N^+(v)\cup \{v\}$, and $N^-[v]$ is defined similarly.
Here is another consequence of \ref{bigandsmall}, more evidence for the truth of \ref{tourconj}:
\begin{thm}\label{bigdom}
For every integer $c\ge 1$, there exists $d\ge 1$ such that if $T$ is a tournament with $\dom(T)\ge d$,
then there is a complete pair $(A,B)$
such that $\chi(A), \chi(B)\ge c$, and $\dom(T[A])\ge c$.
\end{thm}
\Proof
Choose $K,k$ as in \ref{bigandsmall}, and let $d=\max(K,k+c)$.
Let $T$ be a tournament with $\dom(T)\ge d$. Since $d\ge K$, there exists $B\subseteq V(T)$
with $|B|\le k$ such that $\chi(B)\ge c$, by \ref{bigandsmall}. Let $X$ be the set of vertices of $T$ that either belong to $B$ or are
adjacent from some vertex in $B$. Thus $\dom(T\setminus X)\ge c$, since $\dom(T)\ge k+c$
and $|B|\le k$. Let $A=T\setminus X$; then every vertex in $B$ is adjacent from every vertex in $A$, since $A\cap (B\cup X)=\emptyset$. This proves \ref{bigdom}.~\bbox
This reduces proving \ref{tourconj} to proving it for tournaments in which every subtournament has bounded domination number.

The proof of \ref{outnbrs} in~\cite{harut} gives $d$ bounded by a very large, tower-type, function of $c$, but we do not know that it needs to be so
large. In fact we do not have a counterexample for the following, although we find it hard to believe:
\begin{thm}\label{chi/2}
{\bf Conjecture: }For all integers $c\ge 1$ and every tournament $T$ with $\chi(T)\ge 2c$, there exists $v\in V(T)$ such that
$\chi(N^+(v))\ge c$.
\end{thm}
Let us remark that if we replace chromatic number by ``fractional chromatic number'' in the usual sense, then \ref{chi/2} becomes true.
More exactly, let us say a {\em $k$-multicolouring} of $T$ is a list (possibly with repetition) of acyclic subsets of $V(T)$
such that every vertex is in at least $k$ of them, and its {\em size} is the number of sets in the list (counting multiplicity). The
{\em fractional chromatic number} $\chi_f(T)$ is the minimum of $K/k$, over all pairs of positive integers $K,k$ such that there is a $k$-multicolouring of size $K$. Then:
\begin{thm}\label{fractionalchi/2}
For all real numbers $c> 0$ and every tournament $T$ with $\chi_f(T)> 2c$, there exists $v\in V(T)$ such that
$\chi_f(N^+(v))> c$.
\end{thm}
\Proof
Suppose that 
$T[N^+[v]]$ has fractional chromatic number at most $c$, for every vertex $v$. We will prove that $\chi_f(T)\le  2c$. For each $v$, there exist integers 
$K_v, k_v\ge 1$, with $K_v\le  ck_v$, and a $k_v$-multicolouring of $T[N^+[v]]$ of length $K_v$.
Let $K$ be a common multiple of all the numbers $K_v$. For each $v$, by replacing all the sets of the $k_v$-multicolouring of $T[N^+[v]]$ 
by $(K/K_v)k_v$ copies, we obtain a $k_v(K/K_v)$-multicolouring of $T[N^+[v]]$ with size $K$. Thus we may assume that all the 
numbers $K_v$ are equal (to $K$). By replacing all the numbers $k_v$ by $\min(k_v:v\in V(G))$, we may also assume that all the numbers
$k_v$ are equal (to some $k$ say), where $K\le  ck$. Thus, for each $v$, we have a $k$-multicolouring of $T[N^+[v]]$ of size $K$; let
us call it $M_v$. 

By a standard argument using linear programming duality, one can assign a non-negative integer weight $w(v)$ to each vertex $v$ of $T$,
totalling to some positive integer $L$, such that for every vertex $v$, the sum of the weights in $N^-[v]$ is at least $L/2$. 
For each $v\in V(G)$, let us take $w(v)$ copies of $M_v$, and take the union of all these lists, 
forming a list $M$ say. Its size is $KL$. For each vertex $v$, since $\sum_{u\in N^-[v]}w(u)\ge L/2$, and $v$ belongs to $k$ members of $M_u$
whenever $u\in N^-[v]$, it follows that $M$ is a $(kL/2)$-multicolouring of $T$. Since its size is $KL$, it follows that
$\chi_f(T)\le (KL)/(kL/2)=2K/k\le 2c$. This proves \ref{fractionalchi/2}.~\bbox

Two of us proposed~\cite{tourandgraph} a strengthening of \ref{outnbrs}: that if $G$ is a graph with large chromatic number, and $T$ is a tournament with the same vertex set,
then for some vertex $v$, the set of out-neighbours of $v$ in $T$ induces a subgraph of $G$ with large chromatic number. But 
this has very recently been disproved, by 
Gir\~{a}o, Hendrey, Illingworth, Lehner, Michel, Savery and Steiner~\cite{illingworth}, who showed the following:
\begin{thm}\label{tourandgraph}
For all $d\ge 0$, there is a graph $G$ with chromatic number at least $d$, and a tournament $T$ with vertex set $V(G)$, 
such that for each $v\in V(G)$, the set of out-neighbours
of $v$ in $T$ induces a bipartite subgraph of $G$.
\end{thm}

Here is another possible strengthening of \ref{outnbrs} in which we had some hope: that if $A,B$ are disjoint subtournaments of 
a tournament $T$, both with sufficiently large chromatic number, then there is a vertex in one of $A,B$ such that its set of 
out-neighbours in the other set has large chromatic number. But this too is false. 
By a simple modification of their example in \ref{tourandgraph}, Gir\~{a}o et al.~\cite{illingworth} showed:
\begin{thm}\label{bipnbrs}
For all $d\ge 0$, there is a tournament with two disjoint subtournaments $A,B$, both with chromatic number at least $d$, 
such that $\chi(B[N^+(v)])\le 2$ for each vertex $v\in V(A)$, and $\chi(A[N^+(v)])\le 2$ for each vertex $v\in V(B)$.
\end{thm}
Indeed, we suspect that even the following is false, although it remains open for the moment:
\begin{thm}\label{tribipconj}
{\bf Conjecture: }There exists $d\ge 0$ such that if $T$ is a tournament, and $A,B$ are disjoint subsets of $V(T)$
with $\chi(A),\chi(B)\ge d$, then there exist $A'\subseteq A$ and $B'\subseteq B$, both cyclic triangles, such that
one of $(A',B'),(B',A')$ is a complete pair.
\end{thm}

A digression: how much does it matter that we are
concerned with {\em chromatic number} in \ref{chi/2}? If $T$ is a tournament, a {\em submeasure} on $T$ is a function 
$\mu:2^{V(T)}\rightarrow \mathbb{R}^+$ (we use $\mathbb{R}^+$ to denote the set of nonnegative real numbers), 
such that $\mu(\emptyset)=0$, and $\mu$ is increasing and subadditive (that is, $\mu(A)\le \mu(B)$ 
when $A\subseteq B$, and $\mu(A\cup B)\le \mu(A)+\mu(B)$ for all $A,B$). 
Thus chromatic number is a submeasure.
One might hope that we could extend \ref{outnbrs}, to general submeasures instead of chromatic number, but that is false, because of the
following, due to Noga Alon:
\begin{thm}\label{subadd}
For every tournament $T$, there is a submeasure $\mu$ on $T$ such that $\mu(T)$ is the domination number of $T$, and $\mu(N^+(v))=1$
for every vertex $v$.
\end{thm}
\Proof 
For each $X\subseteq V(T)$, let $\mu(X)$ be the cardinality of the smallest subset $Y\subseteq V(T)$ such that every vertex in $X$ either belongs to $Y$ or
is adjacent from a member of $Y$. Then $\mu$ is a submeasure with the desired properties. This proves \ref{subadd}.~\bbox

Let $H$ be a tournament, and for a tournament $T$ define $\chi_H(T)$ to be the minimum $k$ such that
$V(T)$ can be partitioned into $k$ subsets each inducing an $H$-free tournament. 
Thus, $\chi(T)$ is the same as $\chi_H(T)$
when $H$ is a cyclic triangle. Perhaps one can extend \ref{outnbrs} to:
\begin{thm}\label{oneHout}
{\bf Conjecture: }For every tournament $H$ with $\chi(H)\ge 2$ and every integer $c\ge 0$, 
there exists $d\ge 0$ such that for every 
tournament $T$ with $\chi_H(T)\ge d$, there exists $v\in V(T)$ such that
$\chi_H(N^+(v))\ge c$.
\end{thm}

Another digression: 
the domination number of a tournament is at most its chromatic number. One might ask, does \ref{outnbrs} work with chromatic number replaced by domination number? This is false, as a
neat example due to Noga Alon shows:
\begin{thm}\label{onedomout}
For all $d\ge 0$, there is a tournament $T$ with domination number at least $d$, such that $T[N^+(v)]$ has domination number one
for every vertex $v$.
\end{thm}
\Proof
Take a tournament $H$ with domination number at least $d$, and replace each vertex $x$ with a
cyclic triangle $T_x$, such that $T_x\Rightarrow T_y$ if $y$ is adjacent from $x$ in $T$. This forms a larger tournament $T$ say. 
The domination number of $T$ is
at least $d$; but for every vertex $v$ of $T$ , $\dom(T[N^+(v)])=1$ (because if $v\in T_x$, the out-neighbour of $v$ in $T_x$ belongs to and 
dominates $N^+(v)$). This proves \ref{onedomout}.~\bbox

Domination number is very interesting, and we will discuss it more later.
The {\em reverse} of a tournament $T$ is obtained from it by reversing the direction of all its edges.
It is convenient at this point to prove a slight extension of \ref{outnbrs}, the following:
\begin{thm}\label{inandout}
There is an integer-valued function $\phi$ such that for every integer $c\ge 0$, if $T$ is a tournament with $\chi(T)\ge \phi(c)$
then there exists $v\in V(T)$ such that $\chi(N^+(v))\ge c$ and $\chi(N^-(v))\ge c$.
\end{thm}
\Proof
By \ref{outnbrs}, there is an integer-valued function $\psi$ such that for every integer $c\ge 0$, if $T$ is a tournament
with $\chi(T)\ge \psi(c)$,
then there exists $v\in V(T)$ such that $\chi(N^+(v))\ge c$. Define $\phi(c)=2\psi(c)$ for $c\ge 0$. We claim that $\phi$ satisfies
the theorem.

Let $T$ be a tournament with $\chi(T)\ge \phi(c)$, and let $X$ be the set of vertices of $T$ such that $\chi(N^+(v))\ge c$.
It follows that $\chi(T\setminus X)<\psi(c)$,
and so $\chi(X)\ge \chi(T)-\psi(c)\ge \psi(c)$. Hence (by the property of $\psi$, applied to the reverse of 
$T[X]$), there exists $v\in X$ such that $\chi(X\cap N^-(v))\ge c$. But then
$\chi(N^+(v))\ge c$ and $\chi(N^-(v))\ge c$. This proves \ref{inandout}.~\bbox

Let us return to \ref{outnbrs} and weakenings of \ref{tourconj}.
It is a consequence of a result of \cite{heroes} that every tournament with sufficiently large chromatic number
contains a complete pair $(A,B)$ where $A,B$ are both cyclic triangles. 
We will prove a strengthening of this which also is a strengthening of \ref{outnbrs}, and implies \ref{triout}:
\begin{thm}\label{inandouttri}
For all $c\ge 0$ there exists $d\ge 0$ such that if $T$ is a tournament with $\chi(T)\ge d$,
then there exist disjoint sets $P,A,Q\subseteq V(T)$ such that $P\Rightarrow A\Rightarrow Q$, and $A$ is a cyclic triangle, and
$\chi(P),\chi(Q)\ge c$.
\end{thm}
Indeed, this partially extends to $\chi_H$. We will show:
\begin{thm}\label{inandoutH}
For every tournament $H$ with $\chi(H)\ge 2$, and all $c\ge 0$, there exists $d\ge 0$ such that if $T$ is a tournament 
with $\chi_H(T)\ge d$,
then there exist disjoint sets $P,A,Q\subseteq V(T)$ such that $P\Rightarrow A\Rightarrow Q$, where $T[A]$ is isomorphic to $H$, and  
$\chi(P),\chi(Q)\ge c$.
\end{thm}
We have not been able to show the same with the stronger conclusion that $\chi_H(P),\chi_H(Q)\ge c$, though this would be true if conjecture
\ref{oneHout} is true.


\section{Diamonds}\label{sec:diamonds}

A {\em diamond} in a tournament $T$ is a quadruple $(a,b,P,Q)$, where $a,b\in V(T)$ are distinct, and $P,Q$ are
disjoint subsets of $V(T)\setminus \{a,b\}$, such that $a\Rightarrow P\Rightarrow b\Rightarrow Q\Rightarrow a$.
The {\em chromatic number} of a diamond is the minimum of $\chi(P), \chi(Q)$.
Diamonds with large chromatic number are valuable, so in this section
we explore tournaments that
contain no such diamonds.

Fix a numbering $v_1\LL v_n$ of the vertex set of a tournament $T$. For $1\le i,j\le n$, we say $v_j$ is a
{\em right out-neighbour} of $v_i$ if
$j>i$ and $v_j$ is adjacent from $v_i$; and $v_j$ is a {\em right in-neighbour} of $v_i$ if $j>i$ and $v_i$ is adjacent
from $v_j$. Similarly, we say $v_j$ is a {\em left out-neighbour} of $v_i$ if $j<i$ and $v_j$ is adjacent from $v_i$; and
and $v_j$ is a {\em left in-neighbour} of $v_i$ if $j<i$ and $v_i$ is adjacent from $v_j$.
The {\em local chromatic number} of the numbering $v_1\LL v_n$
is the maximum over $1\le i\le n$ of the chromatic number of the set of all $v_j$ such that $v_j$ is either a left
out-neighbour or a right in-neighbour of $v_i$. (Thus, it is the chromatic number of the tournament induced on the set of all neighbours of $v_i$
in the backedge graph under the given numbering. ``Backedge graph'' is defined in section \ref{sec:reduction}.)
Let us show a result that will have several applications, to \ref{heavyedge}, \ref{tworelheroes}, \ref{noStthm}, \ref{newheroes2},
and to the equivalence of \ref{erdosconj} and \ref{tourconj}:

\begin{thm}\label{nodiamond}
If a tournament $T$ admits a numbering with local chromatic number at most $c$, then it contains no diamond with chromatic
number more than $2c$. Conversely, for all $c\ge 0$ there exists $d\ge 0$ such that if a tournament $T$ contains no diamond with
chromatic  number more than $c$, then it admits a numbering with local chromatic number at most $d$.
\end{thm}
\Proof
For the first statement, let $v_1\LL v_n$  be a numbering of $V(T)$ with local chromatic number at most $c$, and let
$(a,b,C,D)$ be a diamond. From the symmetry we may assume that $a=v_i$ and $b=v_j$ where $i<j$ (exchanging $a,b$ and $C,D$
if necessary). If $v_k\in D$, then either $k<j$ and so $v_k$ is a left out-neighbour of $v_j$), or $k>i$
(and then $v_k$ is a right in-neighbour of $v_i$); and so $\chi(D)\le 2c$. This proves the first statement.

For the second, by \ref{outnbrs} there exists $d\ge 0$ such that
if $T$ is a tournament with $\chi(T)> d$, then
there exists $v\in V(T)$ such that $\chi(N^+(v))\ge 2c+2$.
Let $T$ be a tournament that contains no diamond with chromatic number more than $c$. We will show that
$T$ admits a numbering with local chromatic number at most $2d$. Let $H$
be the digraph with vertex set $V(T)$, in which for all distinct $a,b\in V(T)$, $b$ is adjacent from $a$ in $H$ if
$\chi(T[N^+_T(a)\cap N^-_T(b)])\ge 2c+2$.
\\
\\
(1) {\em $H$ has no directed cycle.}
\\
\\
Let $v_1\DD v_2\CC v_k\DD v_1$ be the vertices in order of a directed cycle of $H$.
Since $T$ has no diamond of chromatic number more than $c$ (and therefore, none of order more than $2c$) it follows that
$k>2$. For $1\le i\le k$ let $A_i$ be the set of vertices of $T$ that are adjacent to $v_{i+1}$ and from $v_i$,
where $v_{k+1}$ means $v_1$. Thus $\chi(A_i)\ge 2c+2$ for $1\le i\le k$. For each $i$, let $B_i$
be the set of vertices in $A_i$ adjacent from $v_1$, and $C_i$ be the set of vertices in $A_i$ adjacent to $v_1$. Thus
$B_i\cup C_i = A_i$ if $v_1\notin A_i$, and $B_i\cup C_i = A_i\setminus \{v_1\}$ if $v_1\in A_i$.
Let $I$ be the set of $i\in \{1\LL k\}$ such that $\chi(B_i)>c$, and $J$ the set with $\chi(C_i)>c$. Since
$\chi(B_i\cup C_i)\ge \chi(A_i)-1>2c$, it follows that $i\in I\cup J$ for each $i$. But $1\in I$ and $k\in J$,
and so there exists $i$ with $1\le i<k$ such that $i\in I$ and $i+1\in J$. But then $i(v_1,v_{i+1}, B_i, C_{i+1})$
is a diamond with chromatic number more than $c$, a contradiction. This proves (1).

\bigskip

From (1), there is a numbering $v_1\LL v_n$ of $V(T)$ such that for all $i,j$ with $1\le i<j\le n$,
$v_i$ is not adjacent from $v_j$ in $H$. Let $1\le i\le n$, and let $A$ be the set of right in-neighbours of $v_i$,
and let $B$ be the set of left out-neighbours of $v_i$. For each $v_j\in A$, it follows that $j>i$,
and so $v_jv_i\notin E(H)$, and therefore the set of vertices in $A$ that are out-neighbours of $v_j$
has chromatic number at most $2c+1$. From the choice of $d$, it follows that $\chi(A)\le d$, and similarly
$\chi(B)\le d$ (using in-neighbours instead of out-neighbours). Thus $\{v_1\LL v_n\}$ has chromatic number at most $2d$.
This proves \ref{nodiamond}.~\bbox

For $k\ge 0$, let us say
an edge $uv$ of a tournament $T$ is {\em $k$-heavy} if $\chi(X)\ge k$, where $X$ is the set of all vertices that are adjacent from $v$ and to $u$.
Here is an application of \ref{nodiamond}, proved independently in~\cite{klingel}:
\begin{thm}\label{heavyedge}
For all $k\ge 0$ there exists $q\ge 0$ such that every tournament $T$ with $\chi(T)\ge q$ has a $k$-heavy edge.
\end{thm}
\Proof
By \ref{nodiamond}, there exists $d\ge 0$ such that if a tournament contains no diamond with
chromatic  number at least $k$, then it admits a numbering with local chromatic number at most $d$.
Let $q=2k+4d+1$, and let $T$ be a tournament with $\chi(T)\ge q$.
If $T$ contains a diamond $(a,b,P,Q)$ with chromatic number at least $k$, we may assume that $ab\in E(T)$, by exchanging $a$ with $b$ and $P$ with $Q$
if necessary; but then $ab$ is $k$-heavy. So we assume that there is no such diamond. By \ref{nodiamond}, $T$ admits a numbering $v_1\LL v_n$
with local chromatic number at most $d$. Choose $i_1\le n$ maximum such that $\chi(V_1)\le k+2d$, where $V_1=\{v_1\LL v_{i_1}\}$; if $i_1<n$, choose $i_2$
with $i_1+1\le  i_2\le n$ maximum such that $\chi(V_2)\le k+2d$, where $V_2=\{v_{i_1+1}\LL v_{i_2}\}$, and so on. This partitions $V(T)$ into intervals
$V_1\LL V_t$ say, for some $t\ge 1$; and $\chi(V_s)= k+2d$ for $1\le i<t$, and $\chi(V_t)\le k+2d$. Choose two disjoint sets $P,Q$ both of $k+2d$ colours.
For $1\le s\le t$ with $s$ odd, take a colouring of $T[V_i]$ using colours in $P$, and do the same for $s$ even using colours in $Q$. Since
$\chi(T)>2k+4d$, this is not a valid colouring of $T$, and so there exist $r,s$ with $1\le r,s\le t$ and $s\ge r+2$, such that there is an edge
from $V_s$ to $V_r$. Consequently there exist $i,j$ with $1\le i<j\le n$, such that $v_jv_i\in E(T)$, and $\chi(X)\ge k+2d$,
where $X=\{v_{i+1}\LL v_{j-1}\}$. Let $X_1$ be the set of vertices in $X$ adjacent to $v_i$; let $X_2$ be the set adjacent from $v_j$; and let $X_3$
be the remainder. Since $\chi(X_1), \chi(X_2)\le d$ (because the numbering has local chromatic number at most $d$), it follows that $\chi(X_3)\ge k$,
and so $v_jv_i$ is $k$-heavy. This proves \ref{heavyedge}.~\bbox


\section{Equivalence with the problem of El-Zahar and Erd\H{o}s: the easy half}\label{sec:reduction}

In this section we will show that \ref{tourconj} implies \ref{erdosconj}, and we discuss the converse implication in the next section.
There is a standard technique to go between graphs and tournaments,
as follows.
Let $T$ be a tournament, and choose a numbering $v_1\LL v_n$ of its vertex set. Let $G$ be the graph with vertex set $V(T)$,
in which for $1\le i<j\le n$, the pair $v_i,v_j$ are adjacent in $G$ if and only if $v_i$ is adjacent from $v_j$ in $G$. We call
$G$ the {\em backedge graph} of $T$ under the given numbering. The construction can evidently be reversed: given a graph $G$ and
a numbering, there is a tournament $T$ such that $G$ is the backedge graph of $T$ under the numbering.

We begin with a standard result:
\begin{thm}\label{graphtotour}
Let $G$ be the backedge graph of a tournament $T$ under the numbering $v_1\LL v_n$. Let $\omega(G)$ be the size of the
largest clique of $G$. Then
$$\chi(T)\le \chi(G)\le \omega(G)\chi(T).$$
\end{thm}
\Proof Every set that is stable in $G$ is transitive in $T$, so $\chi(T)\le \chi(G)$. Now let $X$ be transitive in $T$, and let
$<_P$
be the partially ordered set with element set $X$ in which for $v_i,v_j\in X$, we say $v_i<_P v_j$ if $i<j$
and $v_jv_i\in E(T)$. This is indeed a poset, because if $i<j<k$ and $v_i,v_j,v_k\in X$, and
$v_jv_i,v_kv_j\in E(T)$, then $v_kv_i\in E(T)$ (since $X$ is transitive and so $v_iv_k\notin E(T)$). Every totally ordered subset
of the poset is a clique of $G$ and so has size at most $\omega(G)$; and hence by (the dual of) Dilworth's theorem, $X$ can be partitioned
into $\omega(G)$ subsets, each an antichain of the poset and hence each a stable set of $G$. Thus $\chi(G[X])\le \omega(G)$ when $X$
is transitive in $T$; and so $\chi(G)\le \omega(G)\chi(T)$. This proves \ref{graphtotour}.~\bbox

\bigskip

\noindent{\bf Proof of \ref{erdosconj}, assuming \ref{tourconj}.\ \ }
We proceed by induction on $t$; so we may assume that $t\ge 3$, and there exists $d_1$ such that
for every graph $G$ with  $\chi(G)\ge d_1$ and $\omega(G)<t-1$,
there are vertex-disjoint subsets $A,B$ of $V(G)$, both inducing subgraphs with chromatic number at least $c$, with no edges between $A$ and $B$.
Let $c_2=\max(2d_1,2c)$; by the assumed truth of \ref{tourconj}, there is an integer $d_2\ge 1$ such that if $T$ is a tournament
and $\chi(T)\ge d_2$, there is a complete pair $(A,B)$,
where $A,B$ both induce tournaments with chromatic number at least $c_2$.
Let $d=td_2$. Let $G$ be a graph with $\chi(G)\ge d$ and $\omega(G)<t$. We must show that
there are disjoint subsets $A,B$ of $V(G)$, both inducing subgraphs with chromatic number at least $c$, and with no edges beteeen $A,B$.
We may therefore assume that
for every vertex $v\in V(G)$, the subgraph induced on its neighbour set has chromatic number less than $d_1$.

Let $V(G)=\{v_1\LL v_n\}$, and let $T$ be the tournament such that $G$ is its backedge graph under the numbering $v_1\LL v_n$.
From \ref{graphtotour}, it follows that $\chi(T)\ge d/t= d_2$. By \ref{tourconj}, there exist disjoint $A',B'\subseteq V(T)$ with $A'$ complete to $B'$, such that $T[A'],T[B']$
both have chromatic number at least $c_2$.
Choose $i$ minimum such that one of
$$\{v_1\LL v_i\}\cap A', \{v_1\LL v_i\}\cap B'$$
induces a tournament with chromatic number at least $c_2/2$. Now there are two cases.

Suppose first that $\{v_1\LL v_i\}\cap A'$ induces a tournament with chromatic number at least $c_2/2$. Let
$A=\{v_1\LL v_i\}\cap A'$ and $B=\{v_{i+1}\LL v_k\}\cap B'$. Thus $A\ne \emptyset$ and from the minimality of $i$,
$\chi(T\setminus B)\le c_2/2$, and so $\chi(T[B])\ge c_2/2\ge d_1$. Since $A$ is complete to $B$ in $T$,
and $h\le i<j$ for all $h,j$ with $v_h\in A$ and $v_j\in B$, it follows that every vertex of $A$
is adjacent in $G$ to every vertex in $B$. But every subset that is stable in $G$ is acyclic in $T$, and so $\chi(G[B])\ge c_2/2\ge d_1$,
contradicting that for every vertex $v\in V(G)$, the subgraph induced on its neighbour set has chromatic number less than $d_1$.

Thus $\{v_1\LL v_i\}\cap B'$ induces a tournament with chromatic number at least $c_2/2$. Let
$B=\{v_1\LL v_i\}\cap B'$ and $A=\{v_{i+1}\LL v_k\}\cap A'$. As before, $\chi(T[A]), \chi(T[B])\ge c_2/2\ge c$.
Since $A$ is complete to $B$ in $T$,
and $h\le i<j$ for all $h,j$ with $v_j\in A$ and $v_h\in B$, it follows that every vertex of $A$
is nonadjacent in $G$ to every vertex in $B$. Moreover, every subset that is stable in $G$ is acyclic in $T$, and so
$\chi(G[A])\ge \chi(T[A])\ge c$ and similarly $\chi(G[B])\ge c$. This proves \ref{erdosconj}.~\bbox

In the proof above that \ref{tourconj} implies
\ref{erdosconj}, the tournament $T$ constructed
has domination number less than $t$, and so do all its subtournaments. (To see this, choose a clique $X$ of $G$ that is optimal in the following sense:
it contains $v_n$, and it contains $v_j$ where $j<n$ is maximum such that $v_j,v_n$ are adjacent, and it contains $v_i$
where $i<j$ is maximum such that $v_i$ is adjacent to both $v_j,v_n$, and so on. This clique is dominating in $T$.)

So one might consider restricting \ref{tourconj} to tournaments such that all
subtournaments have bounded domination number,
since that would still be strong enough to imply \ref{erdosconj}, as we just saw.

But here is a surprise: that conjecture, apparently much weaker, is the {\em hard} part of \ref{tourconj}: the latter is true for all
tournaments with sufficiently large domination number, because of \ref{bigdom}.

\section{The equivalence of \ref{erdosconj} and \ref{tourconj}: the hard half}\label{sec:eqnce}

As we mentioned earlier, F. Klingelhoefer and A. Newman~\cite{klingel} proved (since an earlier version of this paper was submitted for publication) that the truth of \ref{erdosconj} implies \ref{tourconj}. They used
a result about the chromatic number of digraphs (that is, the smallest $k$ such that the vertex set of the digraph
can be partitioned into $k$ sets each inducing an acyclic digraph), the following:
\begin{thm}\label{klingelthm}
For all integers $a, c\ge 1$ there exists $d$ such that if $D$ is a digraph such that:
\begin{itemize}
\item all directed cycles of $D$ have length at least three;
\item the undirected graph underlying $D$ has no stable set of size $a$; and
\item $D[X]$ has chromatic number at most $c$ for every edge $uv$ of $D$, where $X$ denotes the set of vertices that are adjacent to $u$ and from $v$,
\end{itemize}
then $D$ has chromatic number at most $d$.
\end{thm}
Their derivation of the equivalence of \ref{erdosconj} and \ref{tourconj} from \ref{klingelthm} was quick (about 1.5 pages), but the proof of \ref{klingelthm}
itself was quite lengthy (about 7.5 pages). We have now found a different way to do it, using some ideas from~\cite{klingel} but bypassing \ref{klingelthm},
that we present here.

If $T$ is a tournament, let $J_k$ be the graph with vertex set $V(T)$ and edge set all pairs $\{u,v\}$ such that one of $uv,vu$
is $k$-heavy in $T$. We call $J_k$ the {\em graph of $k$-heavy pairs} of $T$. We begin with:
\begin{thm}\label{lemma1}
Let $T$ be a tournament, and let $d,k\ge 1$ be integers.
If there exists $X\subseteq V(T)$ such that $T[X]$ is transitive and $\chi(J_{2k}[X])\ge  5^d$, then $J_k[X]$ has a $d$-clique.
\end{thm}
\Proof The proof is by induction on $d$, and trivial for $d=1$, so we assume $d\ge 2$ and the result holds for $d-1$. Number $X=\{x_1\LL x_n\}$
such that $x_ix_j\in E(T)$ for $1\le i<j\le n$. (Now the proof is like that of \ref{heavyedge}.) Choose $i_1\le n$ maximum such that $\chi(J_{2k}[X_1])\le 2(5^{d-1})$, where $X_1= \{x_1\LL x_{i_1}\}$; so
either $i_1=n$ or $\chi(J_{2k}[X_1])= 2(5^{d-1})$. If $i_1<n$, choose $i_2$ with $i_1\le i_2\le n$ maximum such that
$\chi(J_{2k}[X_2])\le 2(5^{d-1})$, where $X_2= \{x_{i_1+1}\LL x_{i_2}\}$, and so on. This partitions $X$ into some number of intervals, say $X_1\LL X_t$,
such that $\chi(J_{2k}[X_s])=2(5^{d-1})$ for $1\le s<t$, and $i\chi(J_{2k}[X_t])\le 2(5^{d-1})$. Let us choose two disjoint sets of $2(5^{d-1})$ colours,
say $P,Q$, and colour the graphs $J_{2k}[X_s]$ with colours from $P$ for $1\le s\le t$ with $s$ odd, and use colours from $Q$ to colour the graphs
$J_{2k}[X_s]$  with $s$ even. Since $\chi(J_{2k}[X])\ge 5^d>|P|+|Q|$, this does not make a valid colouring of $J_{2k}[X]$, and so
there exist $1\le r<s\le t$ with $s\ge r+2$ such that there is an edge of $J_{2k}[X]$ between $X_r$ and $X_s$. In particular, there exist
$1\le i\le j\le n$ such that $x_ix_j$ is $2k$-heavy, and $\chi(J_{2k}[X'])\ge 2(5^{d-1})$, where $X'=\{x_{i+1}\LL x_{j-1}\}$.
Let $Y$ be the set of all vertices adjacent from $x_j$ and to $x_i$ in $T$; thus $\chi(Y)\ge 2k$. Let $x\in X'$. Since $\chi(Y)\ge 2k$, either the set of
out-neighbours, or the set of in-neighbours,  of $x$ in $Y$ has chromatic
number at least $k$. In the first case $x_ix$ is $k$-heavy, and in the second $xx_j$ is $k$-heavy. Thus there exist $u\in \{x_i,x_j\}$ and $X''\subseteq X'$
with $\chi(J_{2k}[X''])\ge 5^{d-1}$ such that $u$ is adjacent in $J_k$ to each vertex in $X''$;
and the result follows from the inductive hypothesis applied to $X''$. This proves \ref{lemma1}.~\bbox

In what follows, we will often have a graph ($J_k$ for instance) and a tournament $T$ with the same vertex set, and for $X\subseteq V(T)$, we need to
talk about the chromatic number of the subgraph of $J_k$ induced on $X$, and the chromatic number of the subtournament of $T$ induced on $X$. We denote
the first
by $\chi(J_k[X])$, and the second by $\chi(T[X])$ or just $\chi(X)$.

Let us say a tournament $T$ is {\em $c$-bad} if there do not exist disjoint $A,B\subseteq V(T)$, such that $A\Rightarrow B$ and $\chi(A), \chi(B)\ge c$.
\begin{thm}\label{lemma2}
Let $c,d\ge 1$  be integers.
If \ref{erdosconj} is true, then for all $k\ge c$ there exists $L\ge 0$ such that if $T$ is a $c$-bad tournament with $\chi(T)\ge L$, then its graph $J_k$ of $k$-heavy pairs  has a $d$-clique.
\end{thm}
\Proof
The proof is by induction on $d$, and is trivial for $d=1$, so we assume that $d\ge 2$ and the result holds for $d-1$. Choose $L'$ such that
for every $c$-bad tournament $T$ with $\chi(T)\ge L'$, $J_k$ has a $(d-1)$-clique. Let $K=\max(2k,2L')$. By \ref{nodiamond} there exists $p$ such that
every tournament with no diamond of chromatic number at least $K$ admits
a numbering with local chromatic number at most $p$. By \ref{heavyedge} there exists $q\ge 0$ such that every tournament with chromatic number
at least $q$ has a $k$-heavy edge.

We are assuming the truth of \ref{erdosconj}, and so there exists $M$ such that if a graph $G$ has chromatic number at least $M$ and has no clique of size $d$,
then there are subsets $A,B\subseteq V(G)$ with $\chi(G[A]), \chi(G[B])\ge 5^d(2c+2k+4p)$, such that there are no edges between $A$ and $B$.
Let $L=qM$,
and let $T$ be a $c$-bad tournament with $\chi(T)\ge L$. Suppose first that $T$ contains a diamond $(a,b,P,Q)$ with chromatic number at least $K$.
For each $v\in P$, either the set of out-neighbours, or the set of in-neighbours,  of $v$ in $Q$ has chromatic
number at least $K/2\ge k$; and so one of the edges $av$, $vb$ is $k$-heavy. Consequently, for one of $a,b$ (say $u$), the set of vertices in $P$
adjacent to $u$ in $J_k$ induces a tournament with chromatic number at least $K/2\ge L'$, and so includes a $(d-1)$-clique of $J_k$, from the inductive
hypothesis. Adding $u$
to this clique gives a $d$-clique of $J_k$ as required.

So we assume that there is no such diamond, and consequently $T$ admits a numbering with local chromatic number at most $p$, say $v_1\LL v_n$. We may
assume that $J_k$ has no $d$-clique; and so by \ref{lemma1}, for every subset $X\subseteq V(T)$, $\chi(J_k[X])\le 5^d\chi(X)$.
Conversely, every subset $X\subseteq V(T)$ satisfies $\chi(X)\le q \chi(J_k[X])$, since every set that is stable in $J_k$ induces a tournament
with no $k$-heavy edge and which therefore has chromatic number at most $q$. In particular, $\chi(J_k)\ge \chi(T)/q\ge L/q= M$. By
\ref{erdosconj}, there are subsets $A',B'\subseteq V(T)$ with $\chi(J_k[A']), \chi(J_k[B'])\ge 5^d(2c+2k+4p)$, such that there are no edges of $J_k$
between $A'$ and $B'$. It follows that $\chi(A'), \chi(B')\ge 2c+2k+4p$. Choose $i\le n$ maximum such that $\chi(A_1), \chi(B_1)\le c+k+2p$, where
$A_1=A'\cap \{v_1\LL v_i\}$ and $B_1$ is defined similarly. Let $A_2=A'\setminus A_1$ and $B_2=B'\setminus B_1$. Thus one of $\chi(A_1),\chi(B_1)=c+k+2p$,
and both $\chi(A_2), \chi(B_2)\ge c+k+2p$. We assume $\chi(A_1) = c+k+2p$ (the other case is similar, and indeed the same if we reverse all edges).
Choose $h\le i$ maximum such that $\chi(A)\le c$, where $A=A_1\cap \{v_1\LL v_h\}$; thus $\chi(A)=c$, and $\chi(\{v_{h+1}\LL v_i\})\ge k+2p$. Let $B=B_2$.
As in the proof of \ref{heavyedge}, every edge of $T$ from $B$ to $A$ is $k$-heavy, and so there are no such edges, since
there are no edges of $J_k$
between $A'$ and $B'$. Consequently $A\Rightarrow B$, which is impossible since $T$ is $c$-bad.
This proves \ref{lemma2}.~\bbox

\begin{thm}\label{lemma3}
Let $c\ge 1$. If \ref{erdosconj} is true, then for all integers $d\ge 1$ there exist $L_d, n_d$ such that every $c$-bad tournament with chromatic number at least $L_d$ has a subtournament
with chromatic number at least $d$ and with at most $n_d$ vertices.
\end{thm}
\Proof We use induction on $d$, and can assume that $d>1$, and $L_{d-1}, n_{d-1}$ exist. By \ref{lemma2},
there exists $L_d\ge 0$ such that if $T$ is a $c$-bad tournament with $\chi(T)\ge L_d$, then its graph of $L_{d-1}$-heavy pairs has a $d$-clique.
Define $n_d= d+d(d-1)n_{d-1}/2$. We claim that $L_d,n_d$ satisfy the theorem.
Let $T$ be $c$-bad with $\chi(T)\ge L_d$. By \ref{lemma2}, there is a $d$-clique $X$ of the graph of $L_{d-1}$-heavy pairs. For each edge $uv$ of $T[X]$,
there exists $Y_{uv}\subseteq V(T)$ that is complete to $u$ and from $v$, with $\chi(Y_{uv})\ge L_{d-1}$; and hence from the inductive hypothesis,
there exists $Z_{uv}\subseteq  Y_{uv}$ with $\chi(Z_{uv})\ge d-1$ and $|Z_{uv}|\le n_{d-1}$.
The union of $X$ and all the $Z_{uv}$ has chromatic number at least $d$ (because no two members of $X$ can receive the same colour in a $(d-1)$-colouring),
and has at most $n_d$ vertices. This proves \ref{lemma3}.~\bbox

Finally, we deduce that \ref{erdosconj} implies \ref{tourconj}, because:

\begin{thm}\label{lemma4}
If \ref{erdosconj} is true, then for all $c\ge 1$, every $c$-bad tournament has chromatic number at most  $\max(c2^{n_{2c}}+n_{2c}, L_{2c})$.
\end{thm}
\Proof Suppose that $T$ is a $c$-bad tournament with $\chi(T)\ge  L_{2c})$. By \ref{lemma3}, $T$ has a subtournament $X$
with $\chi(X)\ge 2c$
and with at most $n_{2c}$ vertices. Partition $V(T)\setminus V(X)$ into at most $2^{n_{2c}}$ sets $Y_i\;(i\in I)$ say, such that for each $i\in I$,
all members of $Y_i$ have the same out-neighbours in $V(X)$. Consequently each set $Y_i$ is complete to or from some subset of $V(X)$ with chromatic
number at least $c$, since $\chi(X)\ge 2c$; and so $\chi(Y_i)<c$, since $T$ is $c$-bad. It follows that $\chi(T)\le |X|+2^{n_{2c}}c$. This proves \ref{lemma4}.~\bbox

One can recast the idea of \ref{lemma3} above so that it does not assume the truth of \ref{erdosconj}, as follows:
\begin{thm}\label{newlemma2}
Let $d\ge 1$ be an integer, and let $T$ be a tournament such that no subtournament has chromatic number at least $d$ and has at most $(d!)^2$ vertices. Then for every choice of
integers $L_1\LL L_d\ge 1$, if $\chi(T)\ge L_d$ then there exists $i\in \{2\LL d\}$ and a subset $Y\subseteq V(T)$ such that
$\chi(Y)\ge L_{i}$ and the graph of $L_{i-1}$-heavy edges of $T[Y]$ has no clique of size $i$.
\end{thm}
\Proof
We proceed by induction on $d$. The statement is true when $d=1$, so we assume that $d>1$ and it holds for $d-1$.
We may assume that $\chi(T)\ge L_d$, and 
the graph of $L_{d-1}$-heavy edges of $T$ has a clique of size $d$. So there is a $d$-vertex subtournament $X$ of $T$ such that
all its edges are $L_{d-1}$-heavy. For each edge $uv$ of $X$, let $N(uv)$ be the set of all vertices adjacent to $v$ and from $u$.
If for some edge $uv$ of $X$, there exists $i\in \{2\LL d-1\}$ and a subset $Y\subseteq N_{uv}$ such that
$\chi(Y)\ge L_{i}$ and the graph of $L_{i-1}$-heavy edges of $T[Y]$ has no clique of size $i$, then the theorem is satisfied. Since
$\chi(N(uv))\le L_{d-1}$ for each $uv$, we may therefore assume that for each $uv$, $T[N(uv)]$ has a subtournament $Y_{uv}$
with chromatic number at least $d-1$ 
and at most $((d-1)!)^2$ vertices. Then the union of all the $V(Y_[uv])$ together with $V(X)$ induces a tournament
with chromatic number at least $d$ and with at most $(d!)^2$ vertices. This proves \ref{newlemma2}.~\bbox

It is perhaps of interest to compare this with \ref{noStthm}, which tells us something similar when we exclude one particular 
tournament $\mathcal{S}_t$.

\section{Crossing tournaments}\label{sec:crossing}

Let us say the {\em clique number} of a numbering of a tournament is the size of the
largest clique in its backedge graph.
Can we test whether a tournament admits a numbering with small clique number (even in
an approximate sense) in polynomial time? Is it in co-NP?

In an attempt to show that \ref{erdosconj} implies \ref{tourconj}, before this was proved in~\cite{klingel}, we reduced proving \ref{tourconj} to
proving it for tournaments with small local chromatic number (by applying \ref{noStthm}; we omit the details, since the result is superceded). On the other
hand, \ref{erdosconj} tells us something about tournaments that admit numberings with small clique number.
Could we use it to deduce something about tournaments with small local chromatic number?

Not all numberings with small local chromatic number have small clique number; for instance, the backedge graph
could be a complete graph. But in that case, we could reverse the numbering and get a numbering with small clique number.
That suggests the question, is it true that if a tournament $T$ admits a numbering with small local chromatic number, then it also admits a
numbering with small clique number? The answer is no: we will give a class of tournaments that admit a
numbering with local chromatic number two, such that every numbering has arbitrarily large clique number.

Take a set $P$ of pairs of integers $(a,b)$, such that all the integers
used are distinct (we call this an {\em integer matching}). The graph $H$ with vertex set $P$, in which $(a,b)$ and $(c,d)$
are adjacent if either $a<c<b<d$ or $c<a<d<b$, is called a {\em circle graph}, and these are very interesting graphs;
but we can also derive an interesting tournament 
from $P$, as follows.


Let $T$ be the tournament with vertex set $P$,
in which $(c,d)$ is adjacent from $(a,b)$ if either $a<b<c<d$, or $c<a<b<d$, or $c<a<d<b$. (In other words, if we arrange $a,b,c,d$ in increasing order, then 
the second term is one of $a,b$.) Let us call such a tournament a {\em crossing} tournament. In fact, if we number $P$ by second terms
(so
$(a,b)$ is earlier than $(c,d)$ in the numbering if $b<d$), then the circle graph $H$ is the backedge graph of $T$ under this numbering, 
as we shall see.  
We begin with:
\begin{thm}\label{crossingnumber}
Every crossing tournament admits a numbering such that its backedge graph is a circle graph, and for every vertex, its right in-neighbours are transitive,
and its left out-neighbours are transitive; and consequently the numbering has local chromatic number at most two.
\end{thm}
\Proof
Let $T$ be a crossing tournament defined by  an integer matching $P$. Let
$$P=\{(a_1,b_1)\LL (a_n,b_n)\}$$
where $b_1<b_2<\cdots < b_n$. This defines a numbering of $T$. Let $1\le i,j\le n$. Then under this numbering, $(a_j,b_j)$
is a right in-neighbour of $(a_i,b_i)$ if and only if $a_i<a_j<b_i<b_j$; and so the backedge graph is a circle graph. If
$(a_j,b_j), (a_k,b_k)$ are both right in-neighbours of $(a_i,b_i)$, then the adjacency between them is determined by whichever of
$a_j,a_k$ is larger; $(a_j,b_j)$ is adjacent to $(a_k,b_k)$ if and only if $a_k<a_j$. Consequently
the set of all right in-neighbours of $(a_i,b_i)$ is transitive. Also, $(a_j,b_j)$ is a left out-neighbour of $(a_i,b_i)$
if and only if $a_j<a_i<b_j<b_i$; and so similarly the set of left out-neighbours of $(a_i,b_i)$ is transitive. This proves \ref{crossingnumber}.~\bbox

We do not need the next result; it is included just because we find crossing tournaments interesting.
If $H_1,H_2,H_3$ are tournaments, we denote by $\Delta(H_1,H_2,H_3)$ the tournament $T$ with vertex set the disjoint union of three
sets $A_1,A_2,A_3$, where $T[A_i]$ is isomorphic to $H_i$ for $i = 1,2,3$, and $A_1\Rightarrow A_2\Rightarrow A_3\Rightarrow A_1$.
When $H_3$ is a one-vertex tournament we write $\Delta(H_1,H_2,1)$ for $\Delta(H_1,H_2,H_3)$.
Let $\mathcal{S}_1$ be the tournament with one vertex, and for $t\ge 2$, let
$\mathcal{S}_t=\Delta(\mathcal{S}_{t-1},\mathcal{S}_{t-1},1)$.
\begin{thm}\label{S3crossing}
Crossing tournaments do not contain $\mathcal{S}_3$.
\end{thm}
\Proof
Suppose that $\mathcal{S}_3$ is a crossing tournament, defined by the integer matching $P$ say. Let $P=A\cup B\cup \{c\}$
where $A,B$ are cyclic triangles, and $\{c\}\Rightarrow A\Rightarrow B\Rightarrow\{c\}$.
Let $c=(a,b)$ say. Since $A$ is not transitive, there exists $(p,q)\in A$ such that either $p>a$ or $q<a$; and since
$c$ is adjacent to $(p,q)$, it follows that $p>b$. Similarly there exists $(r,s)\in B$ with $s<a$. But then $(r,s)$
is adjacent to $(p,q)$, contradicting that $A\Rightarrow B$. This proves \ref{S3crossing}.~\bbox

If $P$ is an integer matching, $V(P)$ denotes the set of $2|P|$ ends of its members. If $P,Q$ are integer matchings,
we say that $Q$ is a {\em copy} of $P$ if $|V(P)|=|V(Q)|$ and the (unique) order-preserving bijection from $V(P)$ to $V(Q)$
maps $P$ to $Q$.
We need the following lemma:
\begin{thm}\label{hypergraph}
Let $P$ be an integer matching. Then there is an integer matching $Q$, with the property that for all $Q_1,Q_2$ with $Q_1\cup Q_2=Q$,
there is a copy of $P$ that is a subset of one of $Q_1,Q_2$.
\end{thm}
\Proof Let $|P|=p$ say. We may assume that $V(P)=\{1\LL 2p\}$. Let $N$ be an integer such that for every partition of the
dge set of the complete graph $K_N$ into two classes, there is a $K_{2p}$ subgraph with all edges in the same class.
For $1\le i\le N$ let
$$V_i=\{(i-1)(N-1)+1\LL i(N-1)\};$$
thus, $V_1\LL V_{N}$ form a partition of
$\{1\LL (N-1)N\}$ into sets of cardinality $N-1$. For each pair $(a,b)$ of integers with $1\le a<b\le N$, let
$$f(a,b)=((a-1)(N-1)+b,(b-1)(N-1)+a).$$
Thus the set, $Q$ say, of all these pairs $f(a,b)$ is an integer matching, and for
every choice of $a,b$ with $1\le a<b\le N$, there is a pair $(c,d)\in Q$ with $c\in V_a$ and $d\in V_b$.

We claim that $Q$ satisfies the theorem. Let $Q_1,Q_2\subseteq Q$ with union $Q$. From the choice of $N$, there exists
$I\subseteq \{1\LL N\}$ with $|I|=2p$, such that all the members of $Q$ with both ends in $\bigcup_{i\in I}V_i$
belong to the same one of $Q_1,Q_2$, say to $Q_1$. But then $Q_1$ contains a copy of $P$. This proves \ref{hypergraph}.~\bbox

We deduce:
\begin{thm}\label{cliqueno}
For each integer $k\ge 1$, there is a crossing tournament such that every numbering has clique number at least $k$.
Consequently its chromatic number is at least $k$.
\end{thm}
\Proof
Let us define a crossing tournament $\mathcal{U}_k$ inductively for $k\ge 1$ as follows.
$\mathcal{U}_k$ has one vertex. For $k\ge 2$, we may assume inductively that $\mathcal{U}_{k-1}$ is defined, by the integer matching $P$ say.
By \ref{hypergraph},
there is
an integer matching $Q$ such that for all $Q_1,Q_2$ with union $Q$, one of $Q_1,Q_2$ contains a copy of $P$.

We may assume that $Q$ is a set of pairs of integers
in $\{1\LL 2q\}$. For $i \ge 0$ let $Q^{+i}$ be the integer matching $\{(a+i,b+i):(a,b)\in Q\}$.
Let $b_i=k+1+(2q+1)i$ for $1\le i\le k$, and let
$R= \{(i,b_i):1\le i\le k\}$.
Let $S$ be the union of $R$ and all the sets $Q^{+b_i}$ for $1\le i\le k-1$. Thus $S$ is an integer matching. Let $\mathcal{U}_k$
be the crossing tournament defined by $S$. This completes the inductive definition.

We claim that every numbering of $\mathcal{U}_k$ has clique number at least $k$. The claim is true for $k=1$, so inductively
we assume that $k\ge 2$ and the claim holds for $\mathcal{U}_{k-1}$.
Let $\tau$ be a numbering of $\mathcal{U}_k$,
and let $G$ be the corresponding backedge graph. We need to prove that $T$ has a clique of size $k$.
Let $P,Q,R,S$ and so on be as in the inductive definition
of $\mathcal{U}_k$.

Since $(i,b_i)$ is adjacent from $(j,b_j)$ for $1\le i<j\le k$,
we may assume that there exists $i$ with $1\le i< k$ such that $(i+1,b_{i+1})$ precedes $(i,b_i)$ in the numbering $\tau$, since
otherwise $G$ has a clique of size $k$ as desired.
Every member of $Q^{+b_i}$ is a vertex of $\mathcal{U}_k$, and is a pair $(a,b)$ with $b_i<a<b< b_{i+1}$.
Consequently every member of $Q^{+b_i}$ is either a left out-neighbour of $(i,b_i)$ or a right in-neighbour of $(i+1,b_{i+1})$.
Let $Q_1$ be the set of left out-neighbours of $(i,b_i)$, and $Q_2$ the set of right in-neighbours of $(i+1,b_{i+1})$.
From the choice of $Q$, one of $Q_1,Q_2$, say $Q_j$, contains a copy of $P$; and so there is a clique of $G[Q_j]$
of cardinality $k-1$. If $j=1$ then $(i,b_i)$ is adjacent in $G$ to every vertex of this clique, and if $j=2$
then $(i+1,b_{i+1})$ has the same property. Consequently $G$ contains a clique of size $k$.
This proves that every numbering of $\mathcal{U}_k$ has clique number at least $k$.

If $\mathcal{U}_k$ has chromatic number $t$ say, take a partition into $t$ transitive sets, and take the numbering $\tau$
where we first list the members of the first of the transitive sets, in their natural order, and then list the
members of the second set, and so on. Then the backedge graph is $t$-colourable, and so has no clique of size larger than $t$.
Hence $t\le k$.
This proves \ref{cliqueno}.~\bbox

So admitting a numbering with small local chromatic number does not imply that there is one with
small clique number.


\section{Rebels and posets}\label{sec:heroes}

A tournament $H$ is a
{\em rebel} if for some $c>0$, every $H$-free tournament has domination number less than $c$.
Which tournaments are rebels? 
This section and the next give several new results towards answering this question.

If we arrange the members of a set in a circular order, then every ordered triple of elements rotates clockwise or 
counterclockwise in the natural sense. In particular, if we arrange the vertices of a digraph in a circular order,
then every cyclic triangle is directed clockwise or counterclockwise. 
Let us say a tournament $H$ is a {\em poset tournament} if its vertex set can be arranged in a
circular order, such that there is no clockwise cyclic triangle; or equivalently, if its vertex set
can be numbered $v_1\LL v_n$ such that
for all $i,j,k$ with $1\le i<j<k\le n$,
if $v_iv_j$ and $v_jv_k$ are edges, then $v_iv_k$ is an edge.
Chudnovsky, Kim, Liu, Seymour and Thomass\'e~\cite{tourdom} proved that 
not all tournaments are rebels, and indeed:
\begin{thm}\label{tourdom}
Every rebel is a poset tournament.
\end{thm}
They proposed the conjecture, still open, that
the converse also holds:
\begin{thm}\label{rebel}
{\bf Conjecture: }$H$ is a rebel if and only if $H$ is a poset tournament.
\end{thm}

The conjecture \ref{rebel} is very strong, and here is an entertaining way that one might try to disprove it. Take a class of graphs $\mathcal{F}$, 
that is closed under
taking induced subgraphs, and let 
$\mathcal{F}'$
be the class of all tournaments $T$  that admit a numbering with backedge graph in $\mathcal{F}$. Since $\mathcal{F}'$ is closed under taking subtournaments, if some rebel is not contained in $\mathcal{F}'$  then all members
of $\mathcal{F}'$ have bounded domination number; and if all rebels are contained in $\mathcal{F}'$ then \ref{rebel} would imply that 
every poset tournament is in $\mathcal{F}'$. Consequently
\ref{rebel} would imply that either every poset tournament is in $\mathcal{F}'$, or all members
of $\mathcal{F}'$ have bounded domination number. We tested this on a few familiar classes $\mathcal{F}$, and showed the 
following:
\begin{itemize}
\item If $\mathcal{F}$ is the class of all split graphs (graphs with vertex set the union of a clique and a stable set)
then every member of $\mathcal{F}'$ has domination number at most two.
\item If $\mathcal{F}$ is the class of all line graphs, then every member of $\mathcal{F}'$ has domination number at most three.
\item If $\mathcal{F}$ is the class of all cographs (graphs that do not contain a four-vertex path as an induced subgraph)
then there exists $k$ such that such that all members of $\mathcal{F}'$ have domination number at most $k$.
\item More generally, if $\mathcal{F}$ is the class of all circle graphs (the intersection graph of a set of chords of a circle), there exists $k$
such that all members of $\mathcal{F}'$ have domination number at most $k$.
\item If $\mathcal{F}$ is the class of all graphs with no four-vertex induced cycle,
then every member of $\mathcal{F}'$ has domination number at most four.
\item If $H$ is a permutation graph (that is, there are two numberings $\sigma, \tau$ of $V(H)$, such that $u,v$ are adjacent if $u$
is earlier than $v$ in one of these numberings and later in the other), and $\mathcal{F}$ is the class of all graphs that do not contain $H$ as an induced 
subgraph, then there exists $k$
such that all members of $\mathcal{F}'$ have domination number at most $k$.
\end{itemize}
Here are sketches of the proofs. The first bullet is easy, because tournaments that have a split graph as a backedge graph are two-colourable (because
the subtournaments induced on the clique and on the stable set are both one-colourable).
For the second bullet, the set of in-neighbours of the first vertex in the numbering is the union of two cliques.
The third is a special case of the fourth, because cographs are permutations graphs (easily proved by induction) and hence 
circle graphs. For the fourth, we observe that circle graphs have bounded 
VC-dimension, and hence yield tournaments with bounded VC-dimension, and such tournaments have bounded domination number. 
For the fifth, let $v_1\LL v_n$ be a numbering of a tournament $T$, such that its backedge graph $G$ does not contain 
the four-vertex cycle $C_4$ as an induced subgraph. Choose $1\le i<j\le n$ with $j$ minimum such that $v_j$ is adjacent from $v_i$.
Thus $T[\{v_1\LL v_{j-1}\}]$ is transitive and so has domination number at most one; the set of vertices in $\{v_{j+1}\LL v_n\}$
that are adjacent to both $v_i,v_j$ is transitive, since $G$ does not contain $C_4$; and all other vertices 
in $\{v_{j+1}\LL v_n\}$ are dominated by one of $v_i, v_j$. Hence $\dom(T)\le 4$.
For the sixth, observe that if $H$ is a permutation graph, then it admits a numbering such that the corresponding tournament is transitive, and so the claim follows from \ref{legend} below.

Let us say  a graph $H$ is {\em good} if there exists $k$
such that all members of $\mathcal{F}'$ have domination number at most $k$, where $\mathcal{F}$ is the class of all graphs that do not contain $H$ as 
an induced subgraph.
The last bullet above says that permutation graphs are good, but which other graphs are good?  The complement of any good graph is good,
and the disjoint union of two good graphs is also good (the latter can be shown by an easy modification of the proof above that $C_4$ is good). 
In the other direction, if we set $\mathcal{F}$ to be the class of all 
comparability graphs,
there is no bound on the domination number of the members of $\mathcal{F}'$; and so if $H$ is good, then $\mathcal{F}$ contains $H$, 
and therefore $H$
is a comparability graph. So the class of all good graphs is a superclass of the class of permutation graphs, and a subclass of the class of all comparability graphs; and as far as we know, it might be equal to either one. In particular, here is a cousin of \ref{rebel} (as far as we know, it neither 
implies nor is implied by \ref{rebel}), the following (the ``only if'' part is true):
\begin{thm}\label{rebel2}
{\bf Conjecture: }For every graph $H$, there exists $k$ such that for every tournament $T$ with domination number at least $k$
and every numbering of $T$, the backedge graph contains $H$ as an induced subgraph,
if and only if 
$H$ is a comparability graph.
\end{thm}

Here is another variation, even prettier. We recall that an {\em ordered tournament} is a pair $(T,\tau)$, where $T$ is a tournament
and $\tau$ is a numbering of its vertex set; and its {\em domination number} is the domination number of $T$.
We say an ordered tournament $(H,\sigma)$ 
is a {\em legend} if there exists $k$ such that every ordered tournament $(T,\tau)$ with domination number at 
least $k$ 
contains $(H,\sigma)$ in the natural sense. (Thus a legend is the analogue of a rebel for ordered tournaments.)
Which ordered tournaments 
are legends? We can answer this. 

\begin{thm}\label{legend}
An ordered tournament $(H,\sigma)$ is a legend if and only if $H$ is transitive.
\end{thm}
\Proof
Let us show first that for every legend $(H,\sigma)$, $H$ is transitive. 
The {\em reverse} of an 
ordered tournament is obtained by reversing the direction of all edges, {\em without} reversing the numbering. 
We say a ordered tournament $(H,\sigma)$ is an {\em ordered poset tournament}, where $\sigma$ is the numbering $\sigma_1\LL \sigma_n$,
if $\sigma_k\sigma_i$ is an edge of $H$ for all $1\le i<j<k\le n$, 
such that $\sigma_j\sigma_i$ and $\sigma_k\sigma_j$ are edges of $H$. There are ordered poset tournaments with arbitrarily large domination number,
as was shown in~\cite{tourdom},  
so every legend must be an ordered poset tournament. If we take an ordered poset tournament, and reverse all its edges, and 
reverse the numbering, we obtain another ordered poset tournament; and therefore there are ordered tournaments with 
arbitrarily large domination number, such that  
their reverses are ordered poset tournaments.
Hence, the reverse of every legend must also be an ordered poset tournament. But for
an ordered tournament $(H,\sigma)$, if both $(H,\sigma)$ and its reverse 
are ordered poset tournaments, then $H$ is transitive (as can be seen by checking that no three vertices make a cyclic triangle).
This proves
the ``only if'' part of the theorem.

If $A$ is a subset of $V(T)$, where $T$ is a tournament, we define the {\em external domination number} $\edom(A)$ of $A$ to be the 
size of the smallest subset $X\subseteq V(T)$ such that every vertex in $A\setminus X$ is adjacent from a vertex in $X$.
(This differs from the domination number of the subtournament induced on $A$; and in particular, $\edom$ is subadditive.)
Now let $(H,\sigma)$ be an ordered tournament where $H$ is transitive, and let $\sigma$ be the numbering $\sigma_1\LL \sigma_h$ where $h=|H|$.
We need to show that $(H,\sigma)$ is a legend, but 
for purposes of induction, we will prove a stronger statement, the following:
\\
\\
(1) {\em
If an ordered tournament $(T,\tau)$ has domination number at least $d2^h$ where $d\ge h$ is an integer, and $\tau$ is 
$\tau_1\LL \tau_n$, then there exist $t_1<t_2<\cdots t_h$
such that the restriction of $(T,\tau)$ to $\{\tau_{t_1}\LL \tau_{t_h}\}$ is a copy of $(H,\sigma)$, and each of the $h+1$ sets
$$\{\tau_1\LL \tau_{t_1-1}\},\{\tau_{t_1+1}\LL \tau_{t_2-1}\}\LL \{\tau_{t_h+1}\LL \tau_n\}$$
has external domination number at least $d-h$.}
\\
\\
The claim is true if $h=0$, so we assume that $h\ge 1$ and the result holds for $h-1$. Choose $f$ with $1\le f\le h$ such that
$\sigma_f$ is the vertex of $H$ with in-degree zero; let $H'$ be the transitive tournament 
obtained from $H$ by deleting $\sigma_f$; and let $\sigma'$ be the sequence obtained from $\sigma$
by removing the term $\sigma_f$. In other words, $\sigma'_i=\sigma_i$ for $1\le i<f$, and $\sigma'_i=\sigma_{i+1}$ for $f\le i\le h-1$.
Thus $(H',\sigma')$ is an ordered tournament and $H'$ is transitive.
Let $(T,\tau)$ be an ordered tournament 
with domination number at least $d2^h$, where $d\ge h$ is an integer and $\tau$ is $\tau_1\LL \tau_n$. Temporarily,
define $t_0=0$ and $t_h=n+1$.
From the inductive hypothesis (with $d$ replaced by $2d$)
there exist $1\le t_1<t_2<\cdots t_{h-1}\le n$
such that the restriction of $(T,\tau)$ to $\{\tau_{t_1}\LL \tau_{t_{h-1}}\}$ is a copy of $(H',\sigma')$, and each of the $h$ sets
$\{\tau_{t_{i-1}+1}\LL \tau_{t_{i}-1}\}\; (1\le i\le h)$
has external domination number at least $2d-(h-1)$. 
Let $V=\{\tau_{t_{f-1}+1}\LL \tau_{t_f-1}\}$;
thus, $\edom(V)\ge 2d-(h-1)$. Let $Y$ be the set of
vertices in $V$ that are adjacent from one of the vertices $\tau_{t_1}\LL \tau_{t_{h-1}}$. Since $\edom(Y)\le h-1$, it follows that
$\edom(V\setminus Y)\ge 2d-2(h-1)$. Each vertex $v\in V\setminus Y$ is adjacent to all of $\tau_{t_1}\LL \tau_{t_{h-1}}$, and so
the restriction of  $(T,\tau)$ to $\{\tau_{t_1}\LL \tau_{t_h},v\}$ is a copy of $(H,\sigma)$. But we need to choose $v$ more carefully, 
to arrange the final requirement about external domination number. Since $\edom(V\setminus Y)\ge 2d-2(h-1)$, we may choose 
$k$ minimum such that 
$$\edom(\{\tau_1\LL \tau_{k}\}\cap (V\setminus Y))> d-h.$$ 
The minimality of $k$ implies that $\tau_k\in V\setminus H$, and that 
$$\edom(\{\tau_1\LL \tau_{k-1}\}\cap (V\setminus Y))\le  d-h.$$
Since adding a vertex
to a set changes its external domination number by at most one, it follows that 
$$\edom(\{\tau_1\LL \tau_{k}\}\cap (V\setminus Y))=  d-h+1,$$
and so 
$$\edom(\{\tau_{k+1}\LL \tau_{n}\}\cap (V\setminus Y))\ge (2d-2(h-1))-(d-h+1)= d-h+1.$$
Thus if we insert $k$ into the sequence $t_1\LL t_{h-1}$, so that the new sequence remains increasing, then our requirements are met.
This proves (1).

\bigskip

It follows that if $(T,\tau)$ is a tournament with domination number at least $h2^h$, then it contains $(H,\sigma)$; and
so $(H,\sigma)$ is a legend. This proves \ref{legend}.~\bbox

We remark that if $H$ is transitive and $(H,\sigma)$ is an ordered tournament, then its backedge graph is
a permutation graph, or equivalently the comparability graph of a poset of dimension two. This makes some connection with our 
two ``poset'' conjectures, \ref{rebel} and \ref{rebel2}.

\section{Making rebels}\label{sec:makingrebels}

But, despite all the attacks on it described in the previous section,  \ref{rebel} remains open. 
In support of it, Chudnovsky et al.~\cite{tourdom} proved the following two results:
\begin{thm}\label{twocol}
Every tournament with chromatic number two is a rebel. 
\end{thm}
\begin{thm}\label{D3dom}
$\mathcal{S}_3$ is a rebel.
\end{thm}
The proof of the second was complicated, using a randomized version of a theorem about VC-dimension,  and we will give a much simpler proof.
But first, let us observe:
\begin{thm}\label{chainrebel}
Let $H$ be obtained from the disjoint union of rebels $H_1,H_2$ by making $V(H_1)$ complete to $V(H_2)$. Then $H$ is a rebel.
\end{thm}
\Proof
Choose $d$ such that every tournament with domination number at least $d$ contains both $H_1$ and $H_2$, and let $c=d+|H_2|$. 
Let $T$ be a tournament with $\dom(T)\ge c$; we claim that $T$ contains $H$ and hence $H$ is a rebel.
Since $c\ge d$, there is a copy $S$ of $H_2$ in $T$. Let $X$ be the set of vertices of $T\setminus V(S)$ that are adjacent
to every vertex in $V(S)$. It follows that $\dom(T[X])\ge c-|H_2|$, since every vertex of $T$
not in $X$ either belongs to $V(S)$ or is adjacent from a vertex of $S$. Since $c-|H_2|= d$, we deduce that $T[X]$ contains 
$H_1$, and so $T$ contains $H$. This proves \ref{chainrebel}.~\bbox

A tournament $H$ is a {\em hero} if for some $c>0$, every $H$-free tournament has chromatic number less than $c$.
(Thus,
``hero'' is the concept analogous to ``rebel'' for chromatic number instead of domination number.)
All heroes are rebels, but not all rebels are heroes, and indeed
we know which tournaments are heroes: Berger, Choromanski, Chudnovsky, Fox, Loebl,
Scott, Seymour and Thomass\'e~\cite{heroes} showed that:
\begin{thm}\label{heroes}
A tournament is a hero if and only if all its strongly-connected components are heroes. A strongly-connected
tournament with more than one vertex is a hero if and only if
it equals $\Delta(H,K,1)$ or $\Delta(K,H,1)$ for some hero $H$ and some acyclic tournament $K$.
\end{thm}
The methods of this paper will give a simpler proof of the main part of \ref{heroes}, as we will explain later.

Until now, $\mathcal{S}_3$ was the only non-two-colourable tournament that was known to be a rebel. \ref{chainrebel} trivially 
gives more, but they are not strongly connected.  The next result gives more that are strongly connected, and also
gives a much simpler proof that $\mathcal{S}_3$ is a rebel.
Let us say a {\em ring} in a tournament $T$ is a sequence $X_1,X_2\LL X_n$ of subsets of $V(T)$ with $n\ge 3$, such that
$$X_1\Rightarrow X_2\Rightarrow\cdots\Rightarrow X_n\Rightarrow X_1.$$
(Thus consecutive terms in this sequence are disjoint, but non-consecutive terms may intersect.)
Using rings, we will show that 
if $H_1,H_2,H_3$ are heroes then $\Delta(H_1,H_2, H_3)$ is a rebel. (This implies that $\mathcal{S}_3$ is a rebel, by setting 
$H_1,H_2$ to be cyclic triangles and $H_3$ to have one vertex.)
The claim follows from the following:

\begin{thm}\label{bigtri}
For all integers $c\ge 0$ there exist $K,k$ such that if a tournament $T$ with $\dom(T)\ge K$, then there are three
disjoint sets $A,B,C\subseteq V(T)$ such that $|A|,|B|,|C|\le k$, and $\chi(A), \chi(B),\chi(C)\ge c$, and
$A\Rightarrow B\Rightarrow C\Rightarrow A$.
\end{thm}
\Proof
By \ref{bigandsmall} there exist integers $K,k\ge 1$ such that
every tournament $T$ with $\dom(T)\ge K$ contains a subtournament on at most $k$ vertices having chromatic number
at least $3c$. Also by \ref{bigandsmall} there exist integers $K',k'\ge 1$ such that
every tournament $T$ with $\dom(T)\ge K'$ contains a subtournament on at most $k'$ vertices having chromatic number
at least $2^{2k}c+3c$. We will show that every tournament with domination number at least $\max(K,K')+k+k'$ contains sets $A,B,C$
as in the theorem, each of cardinality at most $k+k'$.

Let $T$ be a tournament with $\dom(T)\ge \max(K',K)+k+k'$.
Let $\mathcal{S}$ be the set of all subsets $X\subseteq V(T)$ with $|X|\le k$ such that $T[X]$ has chromatic number
exactly $3c$.
Let $\mathcal{S}'$ be the set of all subsets $X\subseteq V(T)$ with $|X|\le k'$ such that $T[X]$ has chromatic number
exactly $2^{2k}c+3c$.
Let $\mathcal{R}$ be the set of all subsets that are the union of a member of $\mathcal{S}$ and a member of $\mathcal{S}'$.
For every $X\subseteq V(T)$, if $\dom(T[X])\ge \max(K,K')$, then 
$X$ includes a member of $\mathcal{S}$ and a member of $\mathcal{S}'$, and hence includes a member of $\mathcal{R}$.
In particular, since $\dom(T)\ge \max(K,K')+k+k'$, it follows that $\mathcal{R}\ne \emptyset$.
\\
\\
(1) {\em For each $X\in \mathcal{R}$ there exists $Y\in \mathcal{R}$ with $X\cap Y=\emptyset$ such that $Y\Rightarrow X$.}
\\
\\
Let $X\in \mathcal{R}$. Since $|X|\le k+k'$ and $\dom(T)\ge \max(K',K)+k+k'$, it follows that
the set of vertices in $V(T)\setminus X$ that are complete to $X$ induces a tournament with domination number at least $\max(K,K')$,
and so contains a member of $\mathcal{R}$. 
This proves (1).

\bigskip

From (1) and since $T$ is finite, there is a ring of members of $\mathcal{R}$,
and consequently there is a ring $X_1,X_2\LL X_n$ such that 
\begin{itemize}
\item $|X_1|\le k+k'$ and $\chi(X_1)\ge c$;
\item $|X_2|\le k$ and $\chi(X_2)\ge 2c$;
\item $|X_3|\le k$ and $\chi(X_3)=3c$ (that is, $X_3\in \mathcal{S}$); and
\item $X_i\in \mathcal{R}$ for $4\le i\le n$.
\end{itemize}
Let us call such a sequence a {\em special} ring. Choose a special ring $X_1,X_2\LL X_n$ with $n$ minimum. We may assume that $n\ge 4$, since if $n=3$ then the theorem holds.

Since $X_n$ includes a member of $\mathcal{S}'$ and hence has chromatic number at least $2^{2k}c+3c$,
it follows that $\chi(X_n\setminus X_3)\ge 2^{2k}c$. Moreover, $X_n\cap X_2=\emptyset$, since $X_n\Rightarrow X_1\Rightarrow X_2$.
For each $Z\subseteq X_2\cup X_3$,
let $P_Z$ be the set of vertices in $X_n\setminus X_3$ that are complete to $Z$ and complete from $(X_2\cup X_3)\setminus Z$. 
Since there are
at most $2^{2k}$ choices of $Z$, and each vertex of $X_n\setminus X_3$ belongs to one of the sets $P_Z$, 
it follows that $\chi(P_Z)\ge c$ for some choice of $Z$. If $\chi(X_2\setminus Z)\ge c$, then the theorem holds, since
$$P_Z\Rightarrow X_1\Rightarrow X_2\setminus Z\Rightarrow P_Z;$$
and so we may assume that $\chi(X_2\setminus Z)< c$. Consequently $\chi(X_2\cap Z)\ge c$. If $\chi(X_3\setminus Z)\ge c$
then the theorem holds, since 
$$P_Z\Rightarrow X_2\cap Z\Rightarrow X_3\setminus Z\Rightarrow P_Z;$$
and so we may assume (for a contradiction) that  $\chi(X_3\setminus Z)< c$. Hence $\chi(X_3\cap Z)\ge 2c$ (and consequently
$n\ge 5$, since $X_3$ is not complete to $X_n$). Choose $Y\subseteq X_4$ with $Y\in \mathcal{S}$; then
$$P_Z,X_3\cap Z, Y, X_5\LL X_{n-1}$$
is a special ring, contrary to the minimality of $n$. This proves \ref{bigtri}.~\bbox

We deduce:
\begin{thm}\label{rebellion}
If $H_1,H_2,H_3$ are heroes then $\Delta(H_1,H_2, H_3)$ is a rebel.
\end{thm}
\Proof
Choose $c$ sufficiently large that every tournament with chromatic number at least $c$ contains each of $H_1,H_2,H_3$.
Choose $K,k$ as in \ref{bigtri}. We claim that every tournament $T$ with $\dom(T)\ge K$ contains 
$\Delta(H_1,H_2, H_3)$. Let $T$ be a tournament with $\dom(T)\ge K$. By \ref{bigtri} there exist
$A,B,C$ as in \ref{bigtri}. But $T[A]$ contains $H_1$, and $T[B]$ contains $H_2$, and $T[C]$ contains $H_3$,
and so $T$ contains $\Delta(H_1,H_2, H_3)$. This proves \ref{rebellion}.~\bbox

In view of \ref{chainrebel} and \ref{rebellion}, one might hope for the following:
\begin{thm}\label{trirebel}
{\bf Conjecture: }If $H_1,H_2,H_3$ are rebels then $\Delta(H_1,H_2, H_3)$ is a rebel.
\end{thm}
This is consistent with \ref{rebel}, but we cannot yet prove it.
But here is a special case of \ref{trirebel} that we can prove:
\begin{thm}\label{growrebel}
Let $H$ be a rebel, and let $K$ be a transitive tournament. Then $\Delta(H,K,K)$ is a rebel.
\end{thm}
\Proof
For each integer $r\ge 1$, let $H^r$ be obtained from $r$ disjoint copies of $H$ (say with vertex sets 
$S_1\LL S_r$), by making $V(S_i)$ complete to $V(S_j)$ for $1\le i<j\le r$. If $T[X]$ is isomorphic to $H^r$, we call the 
subsets of $X$ corresponding to $S_1\LL S_r$ the {\em parts} of $X$. By \ref{chainrebel}, $H^r$ is a rebel.
Let $k=|K|$, let $p=k^2$, $q= p+k|H|$, and $r=k(p^k+q^k) |H|^k+1$. 
Choose $c$ such that every tournament with domination number at least $c$ contains $H^r$. Now let $d=c+r|H|$, and let 
$T$ be a tournament
with domination number at least $d$; we will show that $T$ contains $\Delta(H,K,K)$.
Suppose not (for a contradiction).
Let $\mathcal{R}$ be the set of all subsets of $V(T)$ that induce a copy of $H^r$.
Since $d\ge c$, $T$ contains a copy of $H^r$, so $\mathcal{R}\ne \emptyset$; and for each $X\in \mathcal{R}$,
the set of vertices not dominated by $X$ 
(that is, the set of $v\in V(T)\setminus X$ that are adjacent to every vertex of $X$) induces a subtournament
with domination number at least $d-|X|=c$, and consequently also includes a member of $\mathcal{R}$. Hence there is a ring 
$X_1\LL X_n$  of members of $\mathcal{R}$. 

Let us say a {\em special ring} is a ring $X_1\LL X_n$, where $T[X_1]$ is a copy of $H^p$, $T[X_2]$ is a copy of $H$,
and $X_3\LL X_n \in \mathcal{R}$. It follows that $T$ contains a special ring; let us choose a special ring $X_1\LL X_n$
with $n$ minimum. Thus $n\ge 3$. 

Let $P$ be the set of vertices in $X_3$ that have an out-neighbour in at least $k$ parts of $X_1$.
\\
\\
(1) {\em There are fewer than $k(p|H|)^k$ parts of $X_3$ that have a vertex in $P$.}
\\
\\
Suppose that there is a set $Y\subseteq X_3\cap P$ with $|Y|\ge k(p|H|)^k$ such that all vertices in $Y$ belong to different 
parts of $X_3$ (and hence $Y$ is transitive). Each $y\in Y$ has $k$ out-neighbours in $X_1$ that all belong to different 
parts of $X_1$, and hence form a transitive set of cardinality $k$; and there are at most $(p|H|)^k$ choices of such a set, 
since $|X_1|=p|H|$. So there is a subset $Y'\subseteq Y$ with $|Y'|\ge |Y|/(p|H|)^k\ge k$, and a transitive subset
$Z\subseteq X_1$ with cardinality $k$, such that $Y'\Rightarrow Z$. But then $T[X_2\cup Y'\cup Z]$ is a copy of 
$\Delta(H,K,K)$, a contradiction. 
This proves (1).

\bigskip

The set $X_n$ has $r$ parts; let us fix some $q$ of them, called the {\em primary parts} of $X_n$.
Let $Q$ be the set of vertices in $X_3$ that have an out-neighbour in at least $k$ primary parts of $X_n$.
\\
\\
(2) {\em There are fewer than $k(p^k+q^k)|H|^k$ parts of $X_3$ that have a vertex in $P\cup Q$. }
\\
\\
Suppose not; then, by (1), there is a set $Y\subseteq X_3\cap (Q\setminus P)$ with $|Y|\ge k(q|H|)^k$ such that all vertices in $Y$ belong to different
parts of $X_3$ (and hence $Y$ is transitive). As in the proof of (1), there exists $Y'\subseteq Y$ with $|Y'|=k$, and
a transitive subset $Z\subseteq X_n$ with $|Z|=k$, such that $Y'\Rightarrow Z$. Since $Y'\cap P=\emptyset$, each vertex in $Y'$
has an out-neighbour in fewer than $k$ parts of $X_1$, and since $|Y'|=k$ and $X_1$ has $p=k^2$ parts,
there is a part $S$ of $X_1$ such that $S\Rightarrow Y'$. But then $T[S\cup Y'\cup Z]$ is a copy of 
$\Delta(H,K,K)$, a contradiction. 
This proves (2).

\bigskip

Since $r> k(p^k+q^k)|H|^k$, it is not the case that $X_3\Rightarrow X_n$, and so $n\ge 5$. For the same reason, there is a part
$S$ of $X_3$ that is disjoint from $Q$. Each vertex in $S$ has out-neighbours in fewer than $k$ primary parts of $X_n$,
and, since there are $q= k|H|+p$ primary parts of $X_n$, there are $p$ of them that are all complete to $S$. 
Hence there is a copy $T$ of $H^p$ in $T[X_n]$ that is complete to $S$. But then $S, X_4\LL X_{n-1}, T$ is a special ring,
contrary to the minimality of $n$. This proves \ref{growrebel}.~\bbox


\section{Two conjectures of Harutyunyan et al.}\label{sec:harutconj}

Several of our results are based on the
 breakthrough result \ref{bigandsmall} by Harutyunyan et al.~\cite{harut}. In the same paper
they also proposed two strengthenings, that are both still open.
First, they proposed:
\begin{thm}\label{bigandsmalldom}
{\bf Conjecture: }For every integer $c\ge 1$, there exist integers $K,k\ge 1$ such that
every tournament with domination number at least $K$ contains a subtournament on at most $k$ vertices having domination number
at least $c$.
\end{thm}

For the second conjecture, we need some definitions.
We recall that $\mathcal{S}_1$ is the tournament with one vertex, and 
$\mathcal{S}_t=\Delta(\mathcal{S}_{t-1},\mathcal{S}_{t-1},1)$ for $t\ge 2$. For tournaments $T,H$, we say $T$ {\em contains} $H$ if $T$ has a subtournament isomorphic to
$H$.

Harutyunyan et al.~\cite{harut} proposed:
\begin{thm}\label{bigandsmallDk}
{\bf Conjecture: }For every integer $t\ge 1$, there exist $K\ge 1$ such that
every tournament with domination number at least $K$ contains $\mathcal{S}_t$.
\end{thm}

The proof of \ref{bigtri} can be adjusted to show that
conjecture \ref{bigandsmalldom} implies conjecture \ref{bigandsmallDk}. Indeed, if \ref{bigandsmalldom} is true then something stronger than
\ref{bigandsmallDk} holds. Let $\mathcal{T}_1$ be a tournament with one vertex, and for $t>1$ let
$\mathcal{T}_t=\Delta(\mathcal{T}_{t-1},\mathcal{T}_{t-1},\mathcal{T}_{t-1})$. We will show that:

\begin{thm}\label{harutimpl}
Suppose that  for every integer $c\ge 0$, there exists integers $K,k\ge 1$ such that
every tournament with domination number at least $K$ contains a subtournament on at most $k$ vertices having domination number
at least $c$. Then for every integer $t\ge 1$, there exists $K$ such that every $\mathcal{T}_t$-free tournament has domination number at most $K$.
\end{thm}
\Proof
We may assume that $t\ge 2$, and we proceed by induction on $t$. Thus we may assume that there exists $c$  such that every $\mathcal{T}_{t-1}$-free
tournament has
domination number less than $c$.
From the hypothesis, there exist integers $K,k\ge 1$ such that
every tournament $T$ with $\dom(T)\ge K$ contains a subtournament on at most $k$ vertices having domination number
at least $3c$; and there exist integers $K',k'\ge 1$ such that
every tournament $T$ with $\dom(T)\ge K'$ contains a subtournament on at most $k'$ vertices having domination number
at least $2^{2k}c+3c$. We will show that every tournament with domination number at least $\max(K,K')+k+k'$ contains $\mathcal{T}_t$.

Let $T$ be a tournament with $\dom(T)\ge \max(K',K)+k+k'$.
Let $\mathcal{S}$ be the set of all subsets $X\subseteq V(T)$ with $|X|\le k$ such that $\dom(T[X])=3c$; and let
$\mathcal{S}'$ be the set of all subsets $X\subseteq V(T)$ with $|X|\le k'$ such that $\dom(T[X])=2^{2k}c+3c$.

Now the proof proceeds exactly as the proof of \ref{bigtri}: claim (1) holds, and we get a special ring, which now means a ring $X_1,X_2\LL X_n$ such that
\begin{itemize}
\item $|X_1|\le k+k'$ and $\dom(X_1)\ge c$;
\item $|X_2|\le k$ and $\dom(X_2)\ge 2c$;
\item $|X_3|\le k$ and $\dom(X_3)=3c$ (that is, $X_3\in \mathcal{S}$); and
\item $X_i\in \mathcal{R}$ for $4\le i\le n$.
\end{itemize}
We choose such a ring with $n$ minimum, and as
before, we prove that $n=3$. But each of $T[X_1], T[X_2],T[X_3]$ has domination number at least $c$ and so contains $\mathcal{T}_{t-1}$, and hence $T$
contains $\mathcal{T}_t$.
This proves \ref{harutimpl}.~\bbox

It is not easy to think of tournaments with arbitrarily large domination number such that there is some tournament not contained in them.
For instance: 
\begin{itemize}
\item A uniformly random tournament probably has large domination number, but it also probably contains every small tournament.
\item Let $q\ge 3$ be a prime congruent to 3 modulo 4. The {\em Paley tournament} has vertex set the element set
of the field $\mathcal{F}_q$, in which $y$ is adjacent from $x$ if $x-y$ is a square. It was shown by Graham and Spencer~\cite{graham}
that the Paley tournament has domination number at least $\Omega(\log q)$. 
But 
Chung and Graham~\cite{chung} showed that for every tournament $H$, if $q$ is large enough then the Paley tournament contains $H$.
\item Take $2k-1$ linear orderings of the same set $V$, and make a tournament with vertex set $V$
where $v$ is adjacent from $u$ if $v$ is later than $u$ in at least $k$ of the orderings.  This is called a
{\em $k$-majority tournament}. Alon, Brightwell, Kierstead,
Kostochka and Winkler~\cite{kmajority} showed that there are $k$-majority tournaments with domination number at least 
$\Omega(k/ \log k) $. But every tournament is a $k$-majority tournament if $k$ is large enough.
\end{itemize}

Unlike chromatic number, domination number is not monotone; let us
define the {\em subdomination number} of a tournament to be the maximum of
the domination number of all subtournaments. (This is called {\em hereditary domination number} in~\cite{aboulker}.)
The first conjecture \ref{bigandsmalldom} would imply
that we can show that subdomination number is big (when it is) with a constant-time non-deterministic algorithm.
More exactly, for all $k$
there exists $K$, such that if the subdomination number of $T$ is at least $K$, one can demonstrate with a non-deterministic algorithm (with constant
running time if $k$ is fixed)
that its subdomination number is at least $k$. This would be a very nice thing to have, not necessarily via \ref{bigandsmalldom}.

One can show that domination number is small (when it is), just by exhibiting a small dominating set.
But
can we show that subdomination number is small (when it is)?
Is it true that for all
integers $k\ge 0$ there exists $K$, and a poly-time algorithm when $k$ is fixed, that would decide either that $T$
has subdomination number at most $K$ or that $T$ has subdomination number at least $k$? More plausibly, is there a non-deterministic poly-time
algorithm that would do this?

There are many other basic questions about domination number and rebels that we cannot answer; for instance the following three:
\begin{thm}\label{backdom}
{\bf Conjecture: }For all $c\ge 0$, there exists $d\ge 0$ such that if a tournament $T$ has $\dom(T)\ge d$, then it has a subtournament
whose reverse has domination number at least $c$.
\end{thm}
This is trivially true when $c=2$, but we cannot even prove it when $c=3$.
It would imply the following, an analogue of \ref{outnbrs} for subdomination number (the analogue for domination number is false, as we saw in \ref{onedomout}):
\begin{thm}\label{betteronedomout}
{\bf Conjecture: }For all integers $c\ge 0$ there exists $d\ge 0$ such that for every tournament $T$ with subdomination number at least $d$, 
there exists $v\in V(T)$ such that
$T[N^+(v)]$ has subdomination number at least $c$.
\end{thm}
(This conjecture also appears in~\cite{aboulker}.)
\ref{backdom} would also imply:
\begin{thm}\label{rebelreverse}
{\bf Conjecture: }If $H$ is a rebel, then the reverse of $H$ is a rebel.
\end{thm}

\section{The density property.}\label{sec:density}

In this section we prove a result that will be used to deduce \ref{triout} (and some extensions of it), and to prove
\ref{tworelheroes}.
It will be applied to variants of chromatic number, but it holds for general submeasures, and so we have written it in terms of submeasures
(we recall that a submeasure on a tournament $T$ is a function
$\mu:2^{V(T)}\rightarrow \mathbb{R}^+$, such that $\mu(\emptyset)=0$, and $\mu$ is increasing and subadditive).
We call $\mu(X)$ the {\em $\mu$-value} of $X$.

Let $T$ be a tournament, and let $\mu$ be a submeasure on $T$.
If $P,Q\subseteq V(T)$ are disjoint and $c\ge 0$, we denote by 
$\langle P \xrightarrow{c,\mu} Q\rangle$
the set of all $v\in P$
such that
$\mu(N^+(v)\cap Q)\le c$.
Similarly, $ \langle P \xleftarrow{c,\mu} Q\rangle$ denotes the set of all $v\in P$ such that
$\mu(N^-(v)\cap Q)\le c$ (not to be confused with $\langle Q \xrightarrow{c,\mu} P\rangle$). When $\mu$ is the chromatic number $\chi$, we omit the reference to it and write $ \langle P \xrightarrow{c} Q\rangle$, and so on.

Let $g:\mathbb{R}^+\rightarrow\mathbb{R}$ be some function, and let $k \in \mathbb{R}^+$.
A pair $(P,Q)$ of disjoint subsets of $V(T)$ has the {\em $(g,k,\mu)$-out-density property}
if
$$\mu(\langle P \xrightarrow{c,\mu}Q'\rangle)<k$$
for all $c\ge 0$ and for every $Q'\subseteq Q$ with $\mu(Q')\ge g(c)$.

Let $f:\mathbb{R}^+\rightarrow\mathbb{R}^+$ be some function. With $g,k,\mu$
as before, we say a tournament $T$ has the {\em $(f,g,k,\mu)$-branching property} if for all $c\ge 0$, and for every
$X\subseteq V(T)$ with $\mu(X)\ge f(c)$, there is a vertex $a\in X$ and two subsets $P,Q$ of $X$, with the following properties:
\begin{itemize}
\item $P,Q,\{a\}$ are pairwise disjoint, and $Q\Rightarrow \{a\}\Rightarrow P$;
\item $\mu(P),\mu(Q)\ge c$; and
\item $(P,Q)$ has the $(g,k,\mu)$-out-density property.
\end{itemize}
Roughly, this says that in every subset with large $\mu$, we can find a vertex $a$ and two sets $P,Q$ of out-neighbours
and in-neighbours of $a$ respectively, both with large $\mu$, such that
for all $Q'\subseteq Q$ with large $\mu$, there are fewer than $k$ vertices in $P$
with out-neighbour set in $Q'$ of small $\mu$.

How to arrange the $(f,g,k,\mu)$-branching property is a separate issue, and there are combinations of
hypotheses that will give this, that we discuss later. But now we need to prove the following.
\begin{thm}\label{branching}
Let $f,g:\mathbb{R}^+\rightarrow\mathbb{R}^+$ be functions, and let $k\ge 0$. For all $c\ge 0$ and integers $s,t\ge 1$ there
exists $d_{c,s,t}$ with the following property. Let $T$ be a tournament, and let $\mu$ be a submeasure in $T$, such that
$T$ has the $(f,g,k,\mu)$-branching property. Let
$P\subseteq V(T)$, and let $Q_1\LL Q_s$ be subsets of $V(T)\setminus P$ (not necessarily disjoint from each other), with the following properties:
\begin{itemize}
\item $\mu(P), \mu(Q_1)\LL \mu(Q_s)\ge d_{c,s,t}$;
\item for $1\le i\le s$, $(P,Q_i)$ has the $(g,k,\mu)$-out-density property.
\end{itemize}
Then there is a copy $S$ of $\mathcal{S}_t$ in $T[P]$, and for $1\le i\le s$ there is a subset $C_i\subseteq B_i$ complete from $V(S)$
with $\mu(C_i)\ge c$.
\end{thm}
\Proof We prove by induction on $t$ that the statement holds for all choices of $s,c$ and the given value of $t$. Suppose first
that $t=1$. We claim that we may set $d_{c,s,1}=\max(k,g(c))$. Let $T, P, Q_1\LL Q_s$ be as in the theorem. Since $(P,Q_i)$ has
the $(g,k,\mu)$-out-density
property, and $\mu(Q_i)\ge d_{c,s,1}\ge g(c)$, it follows that
$\mu(\langle P \xrightarrow{c,\mu}Q_i\rangle)<k$ for $1\le s$, and so the union of those sets has $\mu$-value
less than $ks$. Since $\mu(P)\ge d_{c,s,1}\ge ks$,
there exists $a\in P$ in none of the sets $\langle P \xrightarrow{c,\mu}Q_i\rangle$, and consequently satisfying the theorem.

Thus the theorem holds when $t=1$; so we assume that $t>1$, and the theorem holds for $t-1$ and all choices of $s,c$.
Let $c'=d_{c,s,t-1}$, and $d' = d_{c',s+1,t-1}$.
We claim that setting $d_{c,s,t}=\max(g(d'), sk+f(d'))$ satisfies the theorem. To see this, let
$T,\mu, P, Q_1\LL Q_s$ be as in the theorem.
For $1\le i\le s$, since $(P,Q_i)$ has the
$(g,k,\mu)$-out-density property, and $\mu(Q_i)\ge d_{c,s,t}\ge g(d')$, it follows that
$$\mu(\langle P \xrightarrow{d',\mu}Q_i\rangle)<k.$$
Consequently the set $P_0$ of vertices in $P$ that belongs to none of the sets $\langle P \xrightarrow{d',\mu}Q_i\rangle$
has $\mu$-value at least $\mu(P)-sk\ge d_{c,s,t}-sk\ge f(d')$.

Since $T$ has the $(f,g,k,\mu)$-branching property, and $\mu(P_0) \ge f(d')$,
there is a vertex $a\in P_0$ and two subsets $P',Q'$ of $P_0$, with the following properties:
\begin{itemize}
\item $P',Q',\{a\}$ are pairwise disjoint, and $Q'\Rightarrow \{a\}\Rightarrow P'$;
\item $\mu(P'),\chi_{\mu}(Q')\ge d'$; and
\item $(P',Q')$ has the $(g,k,\mu)$-out-density property.
\end{itemize}
For $1\le i\le s$, let $Q'_i$ be the set of vertices in $Q_i$ that are adjacent from $a$. So $\mu(Q'_i) \ge d'$ for $1\le i\le s$,
since $a\in P_0$.
Then $Q'$ and $Q_1'\LL Q_s'$ are all disjoint from $P'$; $(P',Q')$ and all the pairs $(P',Q_i')$ have the $(g,k,\mu)$-out-density property;
and $P',Q'$ and all the sets $Q_i'$ have $\mu$-value at least $d'= d_{c',s+1,t-1}$. From the inductive hypothesis,
there is a copy $R$ of $\mathcal{S}_{t-1}$ in $T[P']$, and for $1\le i\le s$ there is a subset $C_i\subseteq Q_i'$ complete from $V(R)$
with $\mu(C_i)\ge c'$, and there is a subset $C\subseteq Q'$ complete from $V(R)$
with $\mu(C)\ge c'$. Now $C_1\LL C_s$ are all disjoint from $C$;
all the pairs $(C,C_i)$ have the $(g,k,\mu)$-out-density property;
and $C$ and all the sets $C_i$ have $\mu$-value at least $c'= d_{c,s,t-1}$.  From the inductive hypothesis,
there is a copy $R'$ of $\mathcal{S}_{t-1}$ in $T[C]$, and for $1\le i\le s$ there is a subset $C_i'\subseteq C_i$ complete from $V(R')$
with $\mu(C_i')\ge c$. But then the subtournament $S$ with vertex set $V(R)\cup V(R')\cup \{a\}$ is a copy of $\mathcal{S}_t$
satisfying the theorem. This proves \ref{branching}.~\bbox

\section{$\MM$-colouring}\label{sec:laws}

Colouring a tournament means partitioning its vertex set into subsets none of which includes a cyclic triangle; and this can
usefully be generalized, as follows. Let $T$ be a tournament, and let $\MM$ be a set of subsets of $V(T)$, each including
the vertex set of a cyclic triangle of $T$. We call $\MM$ a {\em law} for $T$. The {\em order}
of a law is the maximum cardinality of its members (or $|T|$ if $\MM=\emptyset$).
We define the {\em $\MM$-chromatic number}
$\chi_{\MM}(T)$ to be the minimum $k$ such that $V(T)$ can be partitioned into $k$ sets, each including no member of $\MM$.
The function $\chi_{\MM}$ is a submeasure.
Colouring in this sense is the same as colouring the hypergraph $\MM$ in the usual sense of hypergraph colouring;
and $\chi(T)=\chi_{\MM}(T)$ when $\MM$ is the set of all vertex sets of cyclic triangles of $T$.

Now we give the first application of \ref{branching}, to prove \ref{inandoutH} (and therefore its corollaries \ref{inandouttri} and  \ref{triout}),
which we restate:
\begin{thm}\label{Hinandout}
Let $H$ be a tournament and $d\ge 0$. Then there exists an integer $w$ such that for every tournament $T$,
either:
\begin{itemize}
\item there exists a copy $S$ of $H$ in $T$, and two subsets $P,Q\subseteq V(T)\setminus V(S)$ with $Q\Rightarrow V(S)\Rightarrow P$, such that
$\chi(P),\chi(Q)\ge d$; or
\item $V(T)$ can be partitioned into $w$ subsets each inducing an $H$-free subtournament.
\end{itemize}
\end{thm}

This is implied by the following, by taking $\MM$ to be the law of all $X\subseteq V(T)$ such that $T[X]$ is isomorphic to $H$,
and $\ell=|H|$.
\begin{thm}\label{lawinandout}
Let $d,\ell\ge 0$. Then there exists $w$ such that for every tournament $T$, and every law $\MM$ in $T$ of order at most $\ell$,
if  $\chi_{\MM}(T)> w$
then
there exists $A\in \MM$ and two subsets $P,Q\subseteq V(T)\setminus A$ with $Q\Rightarrow A\Rightarrow P$, such that
$\chi(P),\chi(Q)\ge d$.
\end{thm}
(Note that the condition $\chi(P),\chi(Q)\ge d$ refers to chromatic number, and not $\MM$-chromatic number. We do not know if
the same holds using $\MM$-chromatic number.)

We first prove the following weaker statement:
\begin{thm}\label{lawin}
Let $d,\ell\ge 0$. Then there exists $w$ such that for every tournament $T$, and every law $\MM$ in $T$ of order at most $\ell$,
if  $\chi_{\MM}(T)> w$
then
there exist $A\in \MM$ and a subset $Q\subseteq V(T)\setminus A$ with $Q\Rightarrow A$, such that
$\chi(Q)\ge d$.
\end{thm}
\Proof
Let $g$ be the function defined by $g(c)=c\ell +d$ for $c\ge 0$.
By \ref{bigandsmall}, there exist $K,k$ such that every tournament with domination number at least $K$ has a subtournament with chromatic number at least
$d$ and with at most
$k$ vertices.
Let $f$ be the function defined by
$$f(c) = (c+1)\ell+d+ \max(K(c+1),k(c+1)+1)$$
for $c\ge 0$.
With this choice of $f,g$ let $d_{c,s,t}$ be as in \ref{branching}, for all $c,s,t\ge 0$.
Define $w=2d_{2,1,d}$, let $T$ be a tournament, and let $\MM$ be a law in $T$ such that $\chi_{\MM}(T)> w$.
Let $\mu(X)=\chi_{\MM}(X)$ for each $X\subseteq V(T)$; thus $\mu$ is a submeasure.
We suppose for a contradiction that there do not exist $A\in \MM$ and a subset $Q\subseteq V(T)\setminus A$ with $Q\Rightarrow A$, such that
$\chi(Q)\ge d$.
\\
\\
(1) {\em If $P,Q\subseteq V(T)$ are disjoint then $(P,Q)$ has the $(g,2,\mu)$-out-density property.}
\\
\\
Let $c\ge 0$ and let $Q'\subseteq Q$ with $\mu(Q')\ge g(c)$. Suppose that there exists $S\in \MM$ with $S\subseteq
\langle P \xrightarrow{c,\mu}Q'\rangle$. For each $v\in S$, the set of out-neighbours of $v$ in $Q'$ has $\mu$-value at most $c$,
and since $\MM$ has order at most $\ell$, the set of vertices in $Q'$ that have an in-neighbour in $S$ has $\mu$-value at most $c\ell$.
Consequently the set of vertices in $Q'$ that are complete to $S$ has $\mu$-value at least $\mu(Q')-c\ell \ge d$, a contradiction.
Thus there is no such $S$; and so $\mu(\langle P \xrightarrow{c,\mu}Q'\rangle)<2$.
This proves (1).
\\
\\
(2) {\em For every $c\ge 0$ and every $P\subseteq V(T)$ with $\mu(P)\ge (c+1)\ell +d$, there exists $v\in P$ such that $\mu(P\cap N^+(v))\ge c$.}
\\
\\
Choose $S\in \MM$ with $S\subseteq P$. The set of vertices in $P\setminus S$ that are complete to $S$ has chromatic number less than $d$, and hence has
$\MM$-chromatic number less than $d$ (here we use that every member of the law includes a cyclic triangle). So the union over $v\in S$
of the sets $P\cap N^+[v]$ has $\mu$-value at least $\mu(P)-d\ge (c+1)\ell$, and since $|S|\le \ell$, there exists $v\in S$
such that $\mu(P\cap N^+[v])\ge c+1$, and so
$\mu(P\cap N^+(v))\ge c$ (because $\mu$ is subadditive and $\mu(\{v\})\le 1$). This proves (2).
\\
\\
(3) {\em For every $c\ge 0$ and every $P\subseteq V(T)$ with $\mu(P)\ge \max(K(c+1),k(c+1)+1)$, there exists $v\in P$
such that $\mu(P\cap N^-(v))\ge c$.}
\\
\\
From the choice of $K,k$, applied to the reverse of $T[P]$, we deduce that either there is a subset $X\subseteq P$
with $|X|<K$ such that every vertex in $P\setminus X$ has an out-neighbour in $X$, or there exists $C\subseteq P$ with $\chi(C)\ge d$
and $|C|\le k$.
In the first case, since $\mu(P)\ge K(c+1)$, there exists $v\in X$ such that $\mu(P\cap N^-[v])\ge c+1$, and hence
$\mu(P\cap N^-(v))\ge c$
as required. In the second case, the set of vertices in $P\setminus C$ that are complete from $C$ contains no member of $\MM$,
by hypothesis, and so has
$\MM$-chromatic number at most one; and so for some $v\in C$, $\mu(P\cap N^-[v])\ge (\mu(P)-1)/k\ge c+1$, and
hence $\mu(P\cap N^-(v))\ge
c$.
This proves (3).
\\
\\
(4) {\em For every $c\ge 0$ and every $P\subseteq V(T)$ with $\mu(P)\ge f(c)$, there exists $v\in P$ such that
$\mu(P\cap N^+(v))\ge c$ and $\mu(P\cap N^-(v))\ge c$.}
\\
\\
Let $X$ be the set of vertices $v\in P$ such that $\mu(P\cap N^+(v))\ge c$. By (2), $\mu(P\setminus X)< (c+1)\ell +d$, and so
$$\mu(X)\ge \mu(P) -((c+1)\ell+d)\ge \max(K(c+1),k(c+1)+1).$$
Hence the claim follows from (3). This proves (4).

\bigskip

From (4) and (1), it follows that $T$ has the $(f,g,2,\mu)$-branching property.
Since $\mu(T)\ge 2d_{2,1,d}$, there exist disjoint $P,Q\subseteq V(T)$ both with $\MM$-chromatic number at least $d_{2,1,d}$.
By \ref{branching}, with $k=2$, $c=2$, $s=1$ and $t=d$, we deduce that there is a copy $S$ of $\mathcal{S}_d$ in $T[P]$, and there is a subset
$C\subseteq Q$ complete from $V(S)$
with $\mu(C)\ge 2$. Choose $A\in \MM$ with $A\subseteq C$; then $A$ is complete from a copy of $\mathcal{S}_d$, and the latter has chromatic number
at least $d$, a contradiction. This proves \ref{lawin}.~\bbox

We remark that the proof of \ref{lawin} almost never refers to (normal) chromatic number; the only place we need it is when we
apply \ref{bigandsmall}. It would be good to have a version of \ref{lawin} with $\chi(Q)$ replaced by $\chi_{\MM}(Q)$,
but we have not been able to prove this.

Now let us deduce \ref{lawinandout}, which we restate:

\begin{thm}\label{lawinandout2}
Let $d,\ell\ge 0$. Then there exists $w$ such that for every tournament $T$, and every law $\MM$ in $T$ of order at most $\ell$,
if  $\chi_{\MM}(T)> w$
then
there exists $A\in \MM$ and two subsets $P,Q\subseteq V(T)\setminus A$ with $Q\Rightarrow A\Rightarrow P$, such that
$\chi(P),\chi(Q)\ge d$.
\end{thm}
\Proof
Choose $m$ such that for every tournament $T$, and every law $\MM$ in $T$ of order at most $\ell$,
if  $\chi_{\MM}(T)> m$
then
there exists $A\in \MM$ and a subset $P\subseteq V(T)\setminus A$ with $A\Rightarrow P$, such that
$\chi(P)\ge d$. Define $w=m^2$. We will show that $w$ satisfies the theorem.
\\
\\
(1)  {\em For every tournament $T$, and every law $\MM$ in $T$ of order at most $\ell$,
if  $\chi_{\MM}(T)> m$
then
there exists $A\in \MM$ and a subset $P\subseteq V(T)\setminus A$ with $A\Rightarrow P$, such that
$\chi(P)\ge d$.
}
\\
\\
This is immediate from \ref{lawin}, by reversing the direction of all edges.

\bigskip

Now let $T$ be a tournament, and let $\MM$ be a law in $T$ of order at most $\ell$, such that $\chi_{\MM}(T)> w$.
Let $\MM'$ be the set of all $S\in \MM$ such that there exists $P\subseteq V(T)\setminus S$ with $S\Rightarrow P$, and $\chi(P)\ge d$. By (1),
if $X\subseteq V(T)$ and $\chi_{\MM'}(X)\le 1$, then $\chi_{\MM}(X)\le m$; and so, for each $X\subseteq V(T)$,  $\chi_{\MM}(X)\le m\chi_{\MM'}(X)$.
It follows that $\chi_{\MM'}(T)\ge w/m=m$; and hence by \ref{lawin} applied to $\MM'$, there exists $S\in \MM'$ and a subset
$Q\subseteq V(T)\setminus S$ with $Q\Rightarrow
S$,
such that
$\chi(Q)\ge d$. But since $S\in \MM'$ there exists $P\subseteq V(T)\setminus S$ with $S\Rightarrow P$, and $\chi(P)\ge d$. This proves
\ref{lawinandout2}.~\bbox

\section{Excluding $\Delta(H_1,H_2,1)$.}\label{sec:1HH}

We showed in \ref{rebellion} that if $H_1,H_2,H_3$ are heroes then every $\Delta(H_1,H_2, H_3)$-free tournament has bounded 
domination number. But all tournaments that admit numberings with bounded local chromatic number have bounded domination number, as we show below, 
so one might hope to strengthen \ref{rebellion}:
is it true that if $H_1,H_2,H_3$ are heroes then                                       
all $\Delta(H_1,H_2,H_3)$-free tournaments admit numberings 
with bounded local chromatic number? We shall show below that
this is true if one of $H_1,H_2,H_3$ has only one vertex: that is, if $H_1,H_2$ are heroes then all $\Delta(H_1,H_2,1)$-free tournaments admit numberings with bounded local chromatic number.
But otherwise it is false. To see this, observe that if $H_1,H_2,H_3$ each have at least two vertices, 
then $\mathcal{S}_t$ is $\Delta(H_1,H_2,H_3)$-free, and so by \ref{orderTt} of the next section, not all $\Delta(H_1,H_2,H_3)$-free tournaments admit numberings 
with small local chromatic number.

First let us observe that:
\begin{thm}\label{ordertosmalldom}
If a tournament $T$ admits a numbering with local chromatic number at most $c$, then $\dom(T)\le c+1$.
\end{thm}
\Proof
Let $v_1\LL v_n$ be a numbering with local chromatic number at most $c$. Thus $\chi(N^-(v_1))\le c$, and hence
$\dom(T[N^-(v_1)])\le c$ .  But $v_1$ dominates all other vertices of $T$, and so $\dom(T)\le c+1$.
This proves \ref{ordertosmalldom}.~\bbox

A {\em hereditary class} $\mathcal{C}$
is a class of tournaments such that for all $T\in \mathcal{C}$, if a tournament $H$ is isomorphic to a subtournament of $T$, then
$H\in \mathcal{C}$. We can define heroes relative to a hereditary class, as follows. Let $\mathcal{C}$ be a hereditary class, and
let $H$ be a tournament. We say that $H$ is a {\em hero relative to $\mathcal{C}$} if there exists $c\ge 0$ such that
every $H$-free tournament in $\mathcal{C}$
has chromatic number less than $c$. (Consequently, tournaments not in $\mathcal{C}$ are heroes relative to $\mathcal{C}$ if all members of
$\mathcal{C}$
have bounded chromatic number; this is for later technical convenience.) A tournament $H$ is a hero in the earlier sense if and only if it is a hero relative to the class of all tournaments. The greater generality of the current definition will be useful when we explore
$\mathcal{S}_t$-free tournaments in section \ref{sec:St}.

We already defined the  $(g,k,\mu)$-out-density property when $\mu$ is a submeasure; and we say a pair $(P,Q)$ of disjoint subsets 
of $V(T)$ has the
{\em $(g,k,\mu)$-in-density property}
if $\mu(\langle P \xleftarrow{c,\mu}Q'\rangle)<k$ for all $c\ge 0$ and for every $Q'\subseteq Q$ with $\mu(Q')\ge g(c)$.
When $\mu=\chi$ we omit reference to it, and speak of the {\em $(g,k)$-out-density property} and {\em $(g,k)$-in-density
property}, and the
{\em $(f,g,k)$-branching property}.
We will use \ref{branching} to show the next result, which is needed to prove \ref{noStthm}:

\begin{thm}\label{tworelheroes}
Let $\mathcal{C}$ be a hereditary class of tournaments, and let $H_1,H_2$ be heroes relative to $\mathcal{C}$.
There exists $d$ such that for every $\Delta(H_1,H_2,1)$-free tournament $T\in \mathcal{C}$,  every diamond in $T$ has chromatic number
at most $d$. Consequently
there exists $c$ such that every $\Delta(H_1,H_2,1)$-free member of $\mathcal{C}$ admits a numbering
with local chromatic number at most $c$.
\end{thm}
\Proof The first statement implies the second, by \ref{nodiamond}, and so it suffices to prove the first.
Choose $c_0$ such that
every tournament in $\mathcal{C}$ with chromatic number at least $c_0$ contains both $H_1$ and $H_2$. 
Choose $t=c_0+1$.

Let $\phi$ be as in \ref{inandout}. 
Let $g$ be the function defined by $g(c)=c\max(|H_1|,|H_2|)+c_0$ for $c\ge 0$.
Let $f$ be the function defined by $f(c)=\phi(g(c))$ for $c\ge 0$. 
By \ref{branching} (with $\mu=\chi$, and $c=s=0$, and $k=c_0$, and taking $P=V(T)$), there exists $d$
such that every tournament $T$ with $\chi(T)\ge d$ and with the $(f,g,c_0)$-branching property
contains $\mathcal{S}_t$.

Let $T$ be a $\Delta(H_1,H_2,1)$-free member of $\mathcal{C}$. We will show that every diamond in $T$ has chromatic number less than $\max(d,2c_0)$, and
so the theorem holds.
\\
\\
(1) {\em For every diamond $(a,b,P,Q)$ in $T$, the pair $(P,Q)$ has the $(g,c_0)$-out-density property and the 
$(g,c_0)$-in-density property.}
\\
\\
Let $c>0$ and let $Q'\subseteq Q$ with $\chi(Q')\ge g(c)$. We must show that $\chi(\langle P \xrightarrow{c}Q'\rangle)<c_0$
and $\chi(\langle P \xleftarrow{c}Q'\rangle)<c_0$. First we show that $\chi(\langle P \xrightarrow{c}Q'\rangle)<c_0$.

Let $X=\langle P \xrightarrow{c}Q'\rangle$, and suppose that $\chi(X)\ge c_0$. Consequently $T[X]$ includes a copy 
$S$ of $H_2$. For each vertex $v\in V(S)$, the set of vertices in $Q'$ adjacent from $v$ has chromatic number at most $c$, and so the set
of vertices in $Q'$ that have an in-neighbour in $V(S)$ has chromatic number at most $c|H_2|$. Since $\chi(Q')\ge c|H_2|+c_0$,
the set $Y$ of vertices in $Q'$ that are adjacent to each vertex of $S$ has chromatic number at least $c_0$, and so contains
a copy $S'$ of $H_1$. But then the subtournament induced on $V(S)\cup V(S')\cup \{b\}$ is isomorphic to $\Delta(H_1,H_2,1)$,
a contradiction. This proves that $\chi(\langle P \xrightarrow{c}Q'\rangle)<c_0$.

Now we show that $\chi(\langle P \xleftarrow{c}Q'\rangle)<c_0$. (This proof is almost identical to what we just did, 
using $a$ instead of $b$.)
Let 
$X=\langle P \xleftarrow{c}Q'\rangle$, and suppose that $\chi(X)\ge c_0$. Consequently $T[X]$ includes a copy $S$ of $H_1$.
For each vertex $v\in V(S)$, the set of vertices in $Q'$ adjacent to $v$ has chromatic number at most $c$, and so the set
of vertices in $Q'$ that have an out-neighbour in $V(S)$ has chromatic number at most $c|H_1|$. Since $\chi(Q')\ge c|H_1|+c_0$,
the set $Y$ of vertices in $Q'$ that are adjacent from each vertex of $S$ has chromatic number at least $c_0$, and so contains
a copy $S'$ of $H_1$. But then the subtournament induced on $V(S)\cup V(S')\cup \{a\}$ is isomorphic to $\Delta(H_1,H_2,1)$,
a contradiction. This proves that $\chi(\langle P \xleftarrow{c}Q'\rangle)<c_0$, and so proves (1).
\\
\\
(2) {\em
For every diamond $(a,b,P,Q)$ in $T$, if $\chi(Q)\ge 2c_0$ then $T[P]$ has the $(f,g,c_0)$-branching property.}
\\
\\
We must show that for all $c\ge 0$, and for every
$X\subseteq P$ with $\chi(X)\ge f(c)$, there is a vertex $a'\in X$ and two subsets $P',Q'$ of $X$, with the following properties:
\begin{itemize}
\item $P',Q',\{a'\}$ are pairwise disjoint, and $Q'\Rightarrow \{a'\}\Rightarrow P'$;
\item $\chi(P'),\chi(Q')\ge c$; and
\item $(P',Q')$ has the $(g,c_0)$-out-density property.
\end{itemize}
By the definition of $\phi$, and since $\chi(X)\ge f(c)=\phi(g(c))$, there exists $a'\in X$ such that 
$\chi(X\cap N^+(a')), \chi(X\cap N^-(a'))\ge g(c)$. Let $P_1= X\cap N^+(a')$ and $Q_1=X\cap N^-(a')$.
By (1) applied to the diamond $(b,a,Q,P)$, since $\chi(P_1)\ge g(c)$ and $\chi(Q_1)\ge g(c)$, it follows that
$\chi(\langle Q \xleftarrow{c}P_1\rangle)<c_0$ and $\chi(\langle Q \xrightarrow{c}Q_1\rangle)<c_0$.
Since $\chi(Q)\ge 2c_0$, there exists $b'\in Q$ that belongs to neither of the sets
$\langle Q \xleftarrow{c}P_1\rangle, \langle Q \xrightarrow{c}Q_1\rangle$.
Let $P'$ be the set of in-neighbours of $b'$ in $P_1$, and let $Q'$ be the set of out-neighbours of $b'$ in $Q_1$;
thus $\chi(P'),\chi(Q')\ge c$, and $(a',b', P',Q')$ is a diamond. Hence by (1), $(P',Q')$ has the $(g,c_0)$-out-density property,
and so $a',P',Q'$ satisfy the definition of the $(f,g,c_0)$-branching property. This proves (2).

\bigskip

If $\mathcal{S}_{t-1}\in \mathcal{C}$ then  it contains $H_1$ and $H_2$, since its chromatic number is at least $t-1=c_0$; and so
 if $\mathcal{S}_{t}\in \mathcal{C}$ then it contains $\Delta(H_1,H_2,1)$.
Hence 
$T$ does not contain $\mathcal{S}_t$; and it follows from (2) and the choice of $d$ that there is no diamond $(a,b,P,Q)$ in $T$ such that
$\chi(P)\ge d$ and $\chi(Q)\ge 2c_0$. From \ref{nodiamond}, this proves \ref{tworelheroes}.~\bbox


\section{Excluding $\mathcal{S}_t$.}\label{sec:St}

There is no converse to \ref{ordertosmalldom}, because 
the tournament $\mathcal{S}_t$
and all its subtournaments have domination number at most three, and yet:

\begin{thm}\label{orderTt}
Every numbering of $\mathcal{S}_t$ has local chromatic number at least $(t-1)/2$.
\end{thm}
\Proof
We may assume that $t\ge 3$. Let $V(\mathcal{S}_t)$ be the disjoint union of $A,B,\{c\}$, where
$c\Rightarrow A\Rightarrow B\Rightarrow c$
and $A,B$ induce subtournaments isomorphic to $\mathcal{S}_{t-1}$.
Suppose that $v_1\LL v_n$ is a numbering with local chromatic number less than $(t-1)/2$, and choose $i$ minimum
such that one of $A\cap \{v_1\LL v_i\}$,
$B\cap \{v_1\LL v_i\}$ has chromatic number at least $(t-1)/2$. Let $I=\{v_1\LL v_i\}$ and $J=\{v_{i+1}\LL v_n\}$.
Suppose that $\chi(B\cap I)\ge (t-1)/2$. For each $v\in J$,
the set of out-neighbours of $v$ in
$I$ has chromatic number less than $(t-1)/2$ (by assumption, as
they are left out-neighbours), and so $A\cap J=\emptyset$. Hence $\chi(A\cap I)=t-1$, and since $t-2\ge (t-1)/2$,
this contradicts the minimality of $i$. It follows that $\chi(B\cap I)< (t-1)/2$, and so $\chi(A\cap I)\ge (t-1)/2$.
For the same reason it follows
that $c\notin J$, and so $c\in I$. Hence the set of in-neighbours of $c$ in $J$
has chromatic number less than $(t-1)/2$, and so $\chi(B\cap J)<(t-1)/2$. Since $\chi(B\cap I)< (t-1)/2$, we deduce that
$\chi(B)<t-1$, a contradiction. This proves \ref{orderTt}.~\bbox

$\mathcal{S}_t$ is one of the simplest tournaments with large chromatic number, and it would be nice to understand better the tournaments
that do not contain it, particularly since if a tournament contains $\mathcal{S}_t$ with $t$ sufficiently large then it satisfies
\ref{tourconj}. Conjecture \ref{bigandsmallDk} says that all tournaments not containing $\mathcal{S}_t$ have bounded domination number, and for $t=3$,
theorem \ref{tworelheroes} says they all have bounded local chromatic number. 
Does \ref{ordertosmalldom} have a converse if we exclude $\mathcal{S}_t$?
Yes if $t=3$, by \ref{tworelheroes}, but for $t\ge 5$, 
the simplest hope for a converse is false:
for all $t\ge 5$, one can make tournaments that do not contain $\mathcal{S}_t$  that also do not admit numberings with small local
chromatic number. To see this, take a tournament $H$ with large chromatic number that does not contain $\mathcal{S}_3$ 
(for instance, a tournament with backedge graph that has large girth and chromatic number: this has the desired properties);
partition its vertex set into two sets $P,Q$
both with large chromatic number, and add two new vertices $a,b$, where $a\Rightarrow P\Rightarrow b\Rightarrow Q\Rightarrow a$,
forming $T$.
The presence of the diamond guarantees that all numberings of $T$ have large local chromatic number, by \ref{nodiamond}, and yet $T$ 
does not contain $\mathcal{S}_{5}$, since $H$ does not contain $\mathcal{S}_3$.

Here is a more plausible conjecture:
\begin{thm}\label{noStconj}
{\bf Conjecture: }For all integers $t\ge 1$ there exist $K,k$ such that for every $\mathcal{S}_t$-free tournament $T$, there is a partition
of $V(T)$ into $K$ sets each inducing a subtournament that admits a numbering with local chromatic number at most $k$.
\end{thm}
We are far from proving this, because by \ref{nodiamond} it would imply \ref{bigandsmallDk}, but if such $K,k$ exist, then  
in every $\mathcal{S}_t$-free tournament $T$, there is a subtournament with chromatic number at least $\chi(T)/K$ in which every diamond
has chromatic number at most $2k$.
Here is a step in this direction.
The quantifiers are complicated, but it says, roughly, that if $T$ has large chromatic number and does not contain $\mathcal{S}_t$, 
then there is a subtournament with large chromatic number,
in which every diamond has small chromatic number. (We remark that excluding $\mathcal{S}_t$ is necessary: if $T=\mathcal{S}_t$
then every subtournament with large chromatic number has a diamond with large chromatic number.)

\begin{thm}\label{noStthm}
Let $t\ge 2$ be an integer, and let $c_2= 1$. Then 
$$\forall d_2 \exists c_3\forall d_3\exists c_4\cdots \forall d_{t-1}\exists c_t$$
(where $d_2,c_3,d_3\LL c_t$ are all non-negative integers)
such that the following holds. 
If $T$ is an $\mathcal{S}_t$-free tournament, then either every diamond in $T$ has chromatic number less than $c_t$, or
for some $i\in \{2 \LL t-1\}$, there is a subtournament $T'$
with $\chi(T')\ge d_i$,
such that every diamond of $T'$ has chromatic number less than $c_i$.
\end{thm}
\Proof
If $t= 2$, the statement is trivial, since every tournament
with a diamond of positive chromatic number contains $\mathcal{S}_2$. 
Thus, inductively, we may assume that $t\ge 3$ and the result holds for $t-1$. 
Hence 
$$\forall d_2 \exists c_3\forall d_3\exists c_4\cdots \exists c_{t-1},$$
such that if $T$ is a tournament with a diamond of chromatic number at least $c_{t-1}$,
and for each $i\in \{2\LL t-2\}$, every subtournament $T'$
with $\chi(T')\ge d_i$ has 
a diamond of $T'$ with chromatic number at least $c_i$, then $T$ contains $\mathcal{S}_{t-1}$. Let $c_2,d_2\LL c_{t-1}$
be given; then for $d_{t-1}\ge 0$, we need to show that there exists $c_t$ such that the last sentence of the theorem holds.

Let $\mathcal{C}$ be the hereditary class of all tournaments $T$ such that for each $i\in \{2\LL t-1\}$, every subtournament $T'$
with $\chi(T')\ge d_i$ has
a diamond of $T'$ with chromatic number at least $c_i$. Thus 
if $T\in \mathcal{C}$ with $\chi(T)\ge d_{t-1}$, then $T$ has a diamond with chromatic number at least $c_{t-1}$, from the 
definition of $\mathcal{C}$; and from the inductive hypothesis $T$ contains $\mathcal{S}_{t-1}$. Consequently,
$\mathcal{S}_{t-1}$ is a hero relative to $\mathcal{C}_{t-1}$.

From \ref{tworelheroes}, there exists $c_t$ such that for every $\Delta(\mathcal{S}_{t-1},\mathcal{S}_{t-1},1)$-free tournament
$T\in \mathcal{C}$,  every diamond in $T$ has chromatic number
at most $c_t$. Since $\Delta(\mathcal{S}_{t-1},\mathcal{S}_{t-1},1)= \mathcal{S}_{t}$, this proves \ref{noStthm}.~\bbox

\section{Constructing heroes}\label{sec:makingheroes}

Finally, we apply these methods to give a different proof of an old result.
The result of \ref{heroes} gives a construction for all heroes. One half of it, that every hero can be constructed
by the construction given, is easy; the hard half is to show \ref{newheroes1} and \ref{newheroes2} below, and took about eight 
pages in~\cite{heroes}. 
We shall see that, given the machinery developed in this paper, we can easily deduce those two results.

\begin{thm}\label{newheroes1}
If $H$ is obtained from the disjoint union of heroes $H_1,H_2$ by making $V(H_1)$ complete to $V(H_2)$, then $H$ is a hero.
\end{thm}
\Proof
Choose $c_0$ such that every tournament with chromatic number at least $c_0$ contains both $H_1$ and $H_2$. 
Choose $n$ as in \ref{Hinandout}, taking $d=c_0$ and $H=H_1$, and let $c=nc_0$. Now let $T$ be a tournament
with chromatic number at least $c$. 
If $V(T)$ can be partitioned into $n$ subsets each inducing an $H_1$-free subtournament, then each of these subtournaments has 
chromatic number less than $c_0$, and so $\chi(T)< nc_0$, a contradiction.
So, from the choice of $n$, there exists a copy $S$ of $H_1$ in $T$, and a subset $P\subseteq V(T)\setminus V(S)$ with $V(S)\Rightarrow P$, such that
$\chi(P)\ge c_0$. But then $T[P]$ contains $H_2$, and so $T$ contains $H$. This proves \ref{newheroes1}.~\bbox

\begin{thm}\label{newheroes2}
Let $H$ be a hero and $k\ge 0$ an integer: then $\Delta(H,K,1)$ is a hero, where $K$ is a transitive tournament with $k$ vertices.
\end{thm}

This is evidently implied by the following extension of \ref{heavyedge} (since $T[P]$ contains $H$ if $c$ is large enough, and $T[Q]$ contains $K$
if $k$ in \ref{newheroes3} is large enough):
\begin{thm}\label{newheroes3}
For all integers $k,c\ge 0$ there exists $d\ge 0$, such that if $T$ is a tournament with $\chi(T)\ge d$, then 
there exist $v\in V(T)$ and subsets $P\subseteq N^+(v)$ and $Q\subseteq N^-(v)$, where $\chi(P)\ge c$, and $|Q|\ge k$,
and $P\Rightarrow Q$.
\end{thm}
\Proof (Sketch.) Choose $t\ge 2$ such that $\mathcal{S}_t$ contains $v,P,Q$ as specified. 
From \ref{nodiamond} and an easy modification of the
proof of theorem 4.4 of~\cite{heroes}, taking all the sets $X_i$ of that theorem to be singletons, we deduce:
\\
\\
(1) {\em For all $c$, there exists $f(c)$ such that every tournament with chromatic number at least $f(c)$ and with no diamond
of chromatic number at least $c$ contains $v,P,Q$ as specified.}

\bigskip

Now let $T$ be a tournament
that does not contain the desired $v,P,Q$. Consequently it does not contain
$\mathcal{S}_t$, and we must show that its chromatic number is bounded.
By (1), 
for all $c\ge 0$, every subtournament of $T$ with chromatic number at least $f(c)$ has a diamond of chromatic number at least $c$.
From \ref{noStthm}, taking $c_2=1$ and $d_i=f(c_i)$ for $2\le i\le t-1$, we deduce that, with $c_t$ as in that theorem, 
every diamond in $T$ has chromatic number less than $c_t$. Consequently $\chi(T)\le f(c_t)$. This proves \ref{newheroes3}.~\bbox

\section{Acknowledgement}

Our thanks to Noga Alon, Guillaume Aubian and Pierre Aboulker, who gave us counterexamples to some of the conjectures in an earlier version of this paper.

\end{document}